\documentclass[twocolumn]{autart}    
\usepackage[numbers,sort&compress]{natbib}
\usepackage{graphicx,epsfig}
\usepackage{tikz}
\usepackage{amssymb,amsmath,amsfonts}
\usepackage[inline]{enumitem}
\usepackage{epstopdf}
\usepackage{hyperref}

\newcommand{\R}{{\mathbb R}}
\newcommand{\N}{{\mathbb N}}
\newcommand{\Q}{{\mathcal Q	}}
\newcommand{\D}{{\mathcal D	}}
\newcommand{\M}{{\mathcal M	}}
\newcommand{\I}{{\mathcal I	}}
\newcommand{\T}{{\mathcal T	}}
\newcommand{\X}{{\mathcal X	}}
\renewcommand{\S}{{\mathcal S }}
\newcommand{\Vp}{{V_{\text{pm}}}}
\newcommand{\Vpj}{{V_{\text{pm},j}}}
\newcommand{\Vpjj}{{V_{\text{pm},j'}}}
\newcommand{\Vpm}{{V_{\text{pM}}}}
\newcommand{\co}[1]{\overline{co}\left\{#1\right\}}

\newcommand{\fc}{f_{[\nu]}}

\DeclareMathOperator*{\argmax}{arg\,max}
\DeclareMathOperator*{\argmin}{arg\,min}

\begin{document}
\begin{frontmatter}
\title{A Piecewise Smooth Control-Lyapunov Function Framework for Switching Stabilization} 

\author[osu]{Yueyun Lu}\ead{lu.692@osu.edu},
\author[osu]{Wei Zhang}\ead{zhang.491@osu.edu}
\address[osu]{Department of Electrical and Computer Engineering, The Ohio State University, Columbus, OH 43210, USA}
\date{}

\begin{abstract}
This paper studies switching stabilization problems for general switched nonlinear systems. A piecewise smooth control-Lyapunov function (PSCLF) approach is proposed and a constructive way to design a stabilizing switching law is developed. The switching law is constructed via the directional derivatives of the PSCLF with a careful discussion on various technical issues that may occur on the nonsmooth surfaces. Sufficient conditions are derived to ensure stability of the closed-loop Filippov solutions including possible sliding motions. The proposed PSCLF approach contains many existing results as special cases and provides a unified framework to study nonlinear switching stabilization problems with a systematic consideration of sliding motions. Applications of the framework to switched linear systems with quadratic and piecewise quadratic control-Lyapunov functions are discussed and results stronger than the existing methods in the literature are obtained. Application to stabilization of switched nonlinear systems is illustrated through an numerical example. 
\end{abstract}
\end{frontmatter}

\section{Introduction}
This paper studies switching stabilization problems for general switched nonlinear systems in continuous time. The goal is to develop a constructive way to design state-feedback switching laws to ensure closed-loop stabilities including possible sliding motions. The problem is regarded as one of the basic problems in switched systems~\cite{LM99} and has received considerable research attention in recent years~\cite{HML08,LA09,DBPL00,P03,ZAHV09}.

The literature on switched systems has been mainly focusing on their stability analysis~\cite{LiberzonSCL99,liberzon2003switching,ShSIAMR07,LA09,DBPL00}. Depending on the assumptions on switching signal, these studies can be divided into three categories, namely, stability under arbitrary switching, stability under slow switching, and stability under state-dependent switching. Various concepts and tools, such as common Lyapunov functions~\cite{liberzon2003switching}, dwell time~\cite{hespanha1999stability}, and multiple Lyapunov functions~\cite{B98,JR98}, have been proposed to study these problems. These analysis methods typical view the switching inputs as disturbances (e.g. stability under arbitrary or slow switching) or as being generated by a known switching law (e.g. stability under state-dependent switching). They cannot be directly applied to solving stabilization problems, for which one has the freedom to design the switching input to stabilize the system. One main challenge is on the nonsmoothness of switching control laws, and the sliding motions that may occur in closed-loop trajectories. Existing stability analysis results provide little insight into these issues as they typically require the switching signal to be piecewise constant~\cite{liberzon2003switching,DBPL00}, thus eliminating the possibility of having sliding motions that are crucial for switching stabilization problems (see Examples~\ref{ex:unstablesliding} and~\ref{ex:nonuniquesliding} in Section~\ref{sec:formulation}). 

In contrast to stability analysis, switching stabilization has not been adequately studied in the literature. Most of the results focus on quadratic stabilization of switched linear systems (SLS). A switched system is called quadratically stabilizable if it admits a quadratic control-Lyapunov function (CLF)~\cite{ShSIAMR07,SESP99,DBPL00}. It has been shown that a SLS is quadratically stabilizable if there exists a stable convex combination of system matrices~\cite{WPD94,WPD98,LM99,SESP99}. Later, some works extend these quadratic stabilization results by employing piecewise quadratic CLFs. A well-known result along this direction is the largest-region switching strategy proposed in~\cite{LA09,P03}. Such a switching strategy is parameterized by a collection of (possibly overlapping) regions defined in terms of symmetric matrices and a collection of quadratic Lyapunov-like functions. The switching law synthesis problem is then formulated as an optimization problem subject to bilinear matrix inequality (BMI) constraints. The resulting switching law only guarantees closed-loop stability if no sliding motion occurs. A different approach is proposed in a more recent study (\cite{HML08}), for which the switching law is constructed using three classes of composite CLFs defined by taking the pointwise minimum, pointwise maximum, or convex hull of a finite number of quadratic functions. For all of these three classes of composite CLFs, some BMI conditions are derived to ensure closed-loop stability excluding sliding motions. It is further shown that stability including sliding motions can be guaranteed for the case with pointwise minimum CLFs. 

The aforementioned results on switching stabilization have several limitations. Firstly, many of them focus on quadratic stabilization problems, which require the system to have a quadratic CLF. This can be conservative even for SLSs as many SLSs are switching stabilizable without having a quadratic CLF~\cite{liberzon2003switching}. Although several forms of piecewise quadratic CLFs~\cite{P03,HML08,HB98,RJ00} have been considered, the methods are often ad-hoc and rather specific to the particular structure of the adopted piecewise quadratic CLFs. Secondly, many results lack a systematic way to analyze sliding motions. In fact, stability analysis without considering sliding motions may lead to incorrect conclusions about the actual stability behavior of the closed-loop system (see Examples~\ref{ex:unstablesliding} and~\ref{ex:nonuniquesliding} in Section~\ref{sec:formulation}). Lastly, the majority of the existing methods are developed exclusively for SLSs. Switching stabilization of switched nonlinear systems have not been adequately studied. 

In this paper, a piecewise smooth control-Lyapunov function (PSCLF) framework is developed to address the above limitations. We focus on PSCLFs for several reasons. Firstly, they are less restrictive than smooth CLFs in the sense that a switching stabilizable system may not admit a smooth CLF. Secondly, directional derivative always exists for PSCLFs, which is of crucial importance for constructive design of stabilizing switching laws. Lastly, they are general enough to cover a wide range of applications. In fact, most of the existing nonsmooth CLFs proposed in the literature are special classes of PSCLFs. Examples include piecewise quadratic CLFs used in~\cite{WPD94} and the three classes of composite CLFs proposed in~\cite{HML08}. 

In addition to its generality, the proposed PSCLF framework also enables a constructive design of stabilizing switching laws with a systematic consideration of closed-loop sliding motions. In particular, we propose a general procedure to construct a switching law based on a given PSCLF and carefully discuss several technical issues when the sliding surface partly coincides with the nonsmooth surface of the PSCLF. We show that the resulting switching law always guarantees the existence of Filippov solutions of the closed-loop system. We also derive sufficient conditions to ensure closed-loop stability including sliding motions. Different from many existing works, our stability results differentiate sliding motions on attractive sliding surfaces from the ones on nonattractive sliding surfaces. This allows us to focus on sliding motions that are of practical importance, making the overall framework less conservative. The effectiveness of the framework is illustrated through several analytical examples for switched linear systems and a numerical example for a switched nonlinear system. The analytical examples also lead to stronger results than the existing methods on stabilization of SLSs. The proposed PSCLF approach, along with its stability analysis, contains many existing works as special cases and can be used as a unified framework to design stabilizing switching laws with a systematic consideration of sliding motions.

It is worth mentioning that our results on switching stabilization can not be obtained using the classical results on nonsmooth stabilization based on nonsmooth CLFs (e.g. \cite{SS95,CLSS97,clarke1998nonsmooth,BC99,clarke2004lyapunov,S83,sontag1996general,C10,sontag1999stability}). This is due to the fundamental differences between the two classes of problems. (i) Firstly, most classical nonsmooth stabilization results focus on traditional nonlinear systems for which the open-loop vector field $f(x,u)$ is continuous in $(x,u)$. These results are not directly applicable to switched systems whose vector field is not continuous with respect to the switching control. (ii) Secondly, while we focus on Filippov solutions due to its crucial importance for switched systems, many works on nonsmooth stabilization adopt sample-and-hold solutions (e.g. \cite{CLSS97,clarke2004lyapunov,C10,clarke1998nonsmooth}) or Caretheordory solutions (e.g. \cite{AB99}). The difference in the adopted solution notions will significantly change the nature of the problems, making previous methods not directly applicable to our problem. (iii) Lastly, many classical results merely provide sufficient and/or necessary conditions for nonsmooth stabilizability. For example, the classical results on nonsmooth CLF (e.g. \cite{CLSS97,clarke2004lyapunov,C10}) show that asymptotic controllability of a nonlinear system is equivalent to the stabilizability by a measurable feedback law in the sample-and-hold sense. However, there is little indication on how to find such a stabilizing law in a constructive and computationally tractable way. In contrast, this paper offers a simple and constructive way to design the switching law based on a given PSCLF and its corresponding partitions. Due to the above differences, switching stabilization cannot be viewed as a special case of the classical nonsmooth stabilization problems, and the proposed PSCLF framework represents a novel contribution. 

The rest of the paper is organized as follows: In Section 2, we use two examples to illustrate the importance of sliding motions and motivate the problem formulation. In Section 3, we introduce the PSCLF framework and derive two PSCLF theorems to ensure closed-loop stability excluding and including sliding motions, respectively. In addition, we discuss in detail two important special cases with smooth and pointwise minimum CLFs. In Section 4, we use the proposed framework to study stabilization of SLSs under three types of CLFs and provide a numerical example for stabilization of switched nonlinear systems.

\textbf{Notations:} Let $\R_+$ be the set of nonnegative real numbers, $\R^n$ be the n-dimensional Euclidean space. Denoted by $\partial M$ and $\bar{M}$ the boundary and the closure of a set $M \subset \R^n$, respectively. For any positive integer $m$, denoted by $\N_m \triangleq \{1,\cdots,m\}$ the set of positive integers that are less than or equal to $m$. Denoted by $|\cdot|$ the cardinality of a given set, and $\|\cdot\|$ the Euclidean norm of a given vector or matrix. Let $\mu$ be Lebesgue measure. For any $x\in\R^n$ and $\epsilon>0$, define $\mathcal{N}(x;\epsilon)\triangleq \{z\in\R^n: \|z-x\|<r\}$.

\section{Problem Formulation}
\label{sec:formulation}
In this paper, we consider the following continuous-time switched nonlinear system
\begin{align}
\dot{x}(t)=f_{\sigma (t)}(x(t)),\quad \sigma (t) \in \mathcal{Q} \triangleq \{1,\cdots,M\},
\label{eq:ol}
\end{align}
where $x(t)\in \mathbb{R}^n$ denotes the continuous state of the system, $f_i:\mathbb{R}^n \to \mathbb{R}^n$ denotes the vector field of subsystem $i\in \mathcal{Q}$, and $\sigma(t)$ denotes the switching control signal that determines the active subsystem at time $t \in \mathbb{R}_+$. We assume that the vector field of each subsystem $f_i, i \in \mathcal{Q}$ is locally Lipschitz continuous on $\mathbb{R}^n$. In addition, the origin is assumed to be a common equilibrium for all the subsystems $\{f_i\}_{i\in \Q}$.

Assume that the state $x(t)$ is available at all time $t \in \mathbb{R}_+$, and the switching control  is determined through a state-feedback switching law $\nu:\mathbb{R}^n \to \mathcal{Q}$. The corresponding closed-loop system can be written as
\begin{align}
\label{eq:cl}
\dot{x}(t)=f_{\nu(x(t))}(x(t)).
\end{align}
For simplicity, the closed-loop vector field under a switching law $\nu$ will be denoted  by $\fc$, i.e., $\fc(x)\triangleq f_{\nu(x)}(x)$, $\forall x\in\R^n$. Although each subsystem vector field $f_i$ is continuous, a switching law $\nu$ may introduce discontinuities in the closed-loop vector field $\fc$. In general, the differential equation in (\ref{eq:cl}) may not have a classical or Caratheodory solution \cite{C09}. In this paper, we adopt the Filippov solution notion \cite{filippov1988differential} to study the switching stabilization problem.

\begin{defn} \label{def:filippovmap} \textit{(Filippov Set-Valued Map \cite{C09})} For any vector field $X: \mathbb{R}^n \to \mathbb{R}^n$, the corresponding Filippov set-valued map $F[X]:\mathbb{R}^n \to \mathfrak{B}(\mathbb{R}^n)$ is defined as 
\begin{align*}
F[X](x)\triangleq \bigcap_{\delta>0} \bigcap_{\mu(S)=0} \co{X(\mathcal{N}(x;\delta)\backslash S)}, \quad x \in \mathbb{R}^n,
\end{align*}
where $\mathfrak{B}(\mathbb{R}^n)$ denotes the collection of subsets of $\mathbb{R}^n$, $\overline{co}$ denotes convex closure and $\mu$ denotes Lebesgue measure. 
\end{defn}

\begin{defn} \textit{(Filippov Solution \cite{C09})} A Filippov solution to a differential equation $\dot{x}(t)=X(x(t))$ over $[0,t_1]$ with $t_1>0$ is an absolutely continuous map $x:[0,t_1]\to \mathbb{R}^n$ that satisfies the differential inclusion $\dot{x}(t)\in F[X](x(t))$ for almost all $t \in [0,t_1]$.
\end{defn}

According to Definition \ref{def:filippovmap}, we can exclude an arbitrary set of measure zero around $x$ when computing $F[X](x)$. This implies that two vector fields that differ on a set of measure zero will lead to the same Filippov set-valued map, and hence the same Filippov solution. In addition, it can be verified that if the vector field $X$ is continuous, then the Filippov solution to $\dot{x}(t)=X(x(t))$ coincides with the classical solution. If $X$ is discontinuous, then the Filippov solution exists as long as the map $X:\mathbb{R}^n \to \mathbb{R}^n$ is measurable and locally essentially bounded \cite{filippov1988differential}. Since we assume that the vector field of each subsystem is continuous, it can be easily verified that the Filippov solution to the closed-loop system (\ref{eq:cl}) exists whenever the switching law $\nu:\mathbb{R}^n \to \mathcal{Q}$ is measurable. However, switching laws may take unappealing forms if we merely require them to be measurable. For the purpose of systematic analysis on sliding motions, we focus on a class of switching laws defined as follows.

\begin{defn} \label{def:switchinglaw} \textit{(Admissible Switching Law)} A switching law $\nu:\R^n\to\Q$ is called \emph{admissible} if there exists a collection of disjoint and open sets $\{\D_i\}_{i\in\Q}\subseteq \R^n$ satisfying $\cup_{i\in\Q}\bar{\D}_i=\R^n$ such that $\nu(x)=i$, $\forall x\in\D_i$.
\end{defn} 

The sets $\{\D_i\}_{i\in\Q}$ and $\{\partial\D_i\}_{i\in\Q}$ in Definition \ref{def:switchinglaw} will be referred to as the {\em switching partitions} and the {\em switching boundaries}, respectively. Let $\partial \D \triangleq \cup_{i\in\Q}\partial\D_i$. In the above definition, each switching partition $\D_i$ is not required to be connected. In addition, since $\{\D_i\}_{i\in\Q}$ are full-dimensional open subsets of $\R^n$, the switching boundaries $\{\partial \D_i\}_{i\in\Q}$ are of measure zero. As a result, an admissible switching law $\nu$ is piecewise constant. It is easy to see that an admissible switching law ensures the closed-loop vector field $\fc$ to be measurable and locally essentially bounded. Thus, a Filippov solution to system~(\ref{eq:cl}) always exists. 

We denote $x(\cdot;z,\nu):\mathbb{R}_+ \to \mathbb{R}^n$ as a Filippov solution to the closed-loop system (\ref{eq:cl}) under an admissible switching law $\nu$ with initial state $z \in \mathbb{R}^n$. Sliding motion may occur along a Filippov solution \cite{utkin1992sliding}. More precisely, if there exists a nontrivial time interval $(t_1,t_2)$ such that $\dot{x}(t;z,\nu) \neq f_{[\nu]}(x(t;z,\nu))$ for almost all $t \in (t_1,t_2)$, then this part of the state trajectory is called a sliding motion where the velocity of the state $\dot{x}(t;z,\nu)$ can take any value within the convex set $F[f_{[\nu]}](x(t;z,\nu))$.

Since the closed-loop vector field $\fc$ is continuous inside each of the switching partition $\D_i$, a sliding motion of $\fc$ can only occur on the switching boundaries $\partial\D$. For each $x\in \partial\D$, let $\I_{sm}(x)$ be the set of indices of all the subsystems that may possibly be involved in a sliding motion, i.e., 
\begin{align}
\I_{sm}(x) \triangleq \{i\in\Q: x\in\partial\D_i\}.
\label{eq:Ism}
\end{align}
During sliding motion, the velocity satisfies $\dot x = \sum_{i\in\I_{sm}(x)} \alpha_i(x) f_i(x)$ for some convex combination coefficients $\{\alpha_i(x)\}_{i\in\I_{sm}(x)}$ with $\sum_{i\in\I_{sm}(x)}\alpha_i(x) =1$. These coefficients can be determined if the sliding surface equation is known~\cite{DBPL00, U77}. In addition, we define $\I_{sm}^a(x)$ as the set of indices of all the \textit{active} subsystems involved in the sliding motion, i.e., $\I_{sm}^a(x)\triangleq \{i\in \I_{sm}(x):\alpha_i(x)\neq 0 \text{ where } \dot{x}=\sum_{i\in\I_{sm}(x)}\alpha_i(x)f_i(x), \sum_{i\in\I_{sm}(x)} \alpha_i(x)=1\}$. Obviously, $\I_{sm}^a(x)\subseteq\I_{sm}(x)$, $\forall x\in \partial\D$.

Existing works in the literature often discuss stability properties by excluding sliding motions. However, this may lead to an incorrect conclusion about the actual stability behavior as illustrated in the following example.

\begin{figure*}[t!]
	\centering
	\begin{tabular}{cc}
		\includegraphics[width=0.35\textwidth]{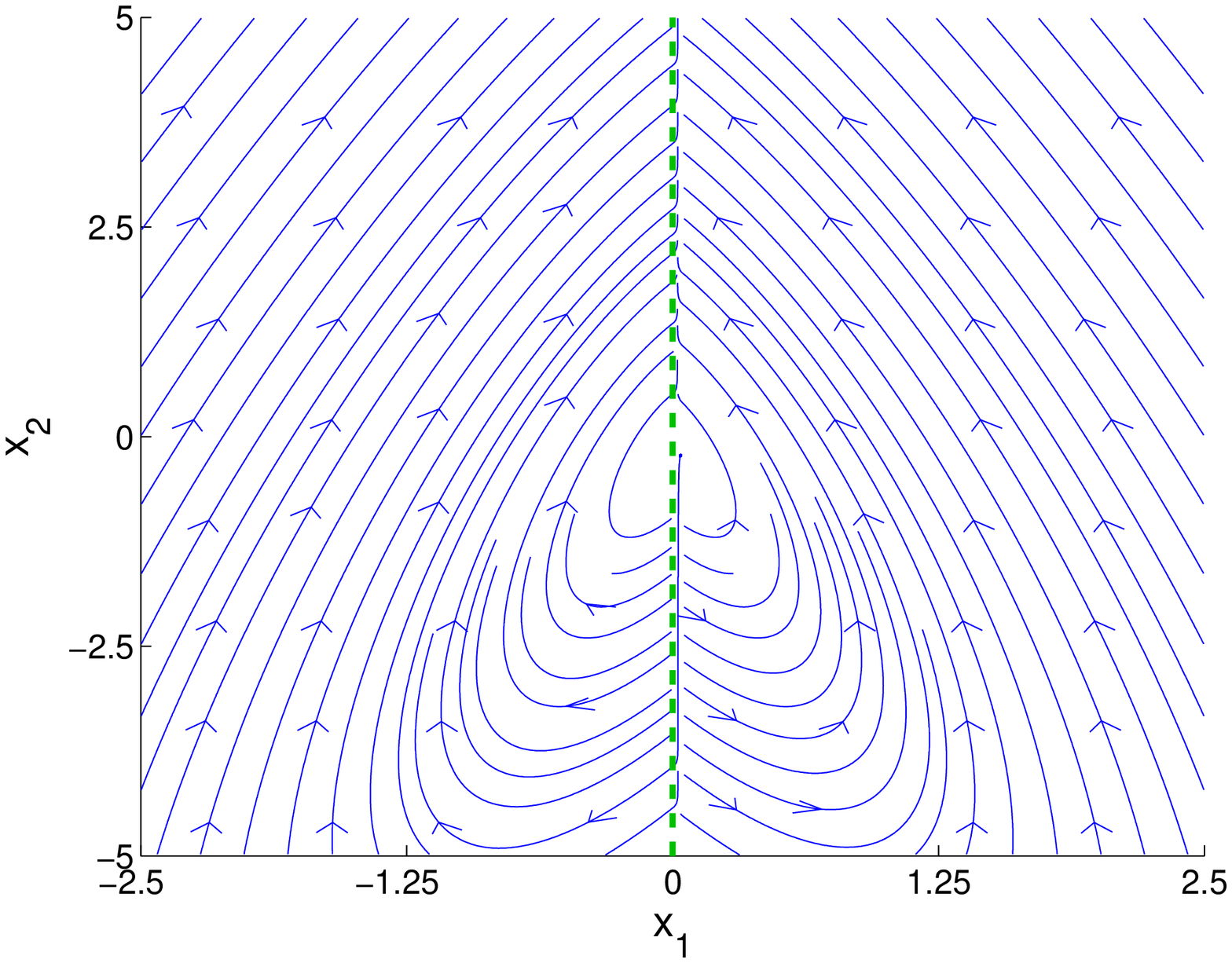} &
		\includegraphics[width=0.35\textwidth]{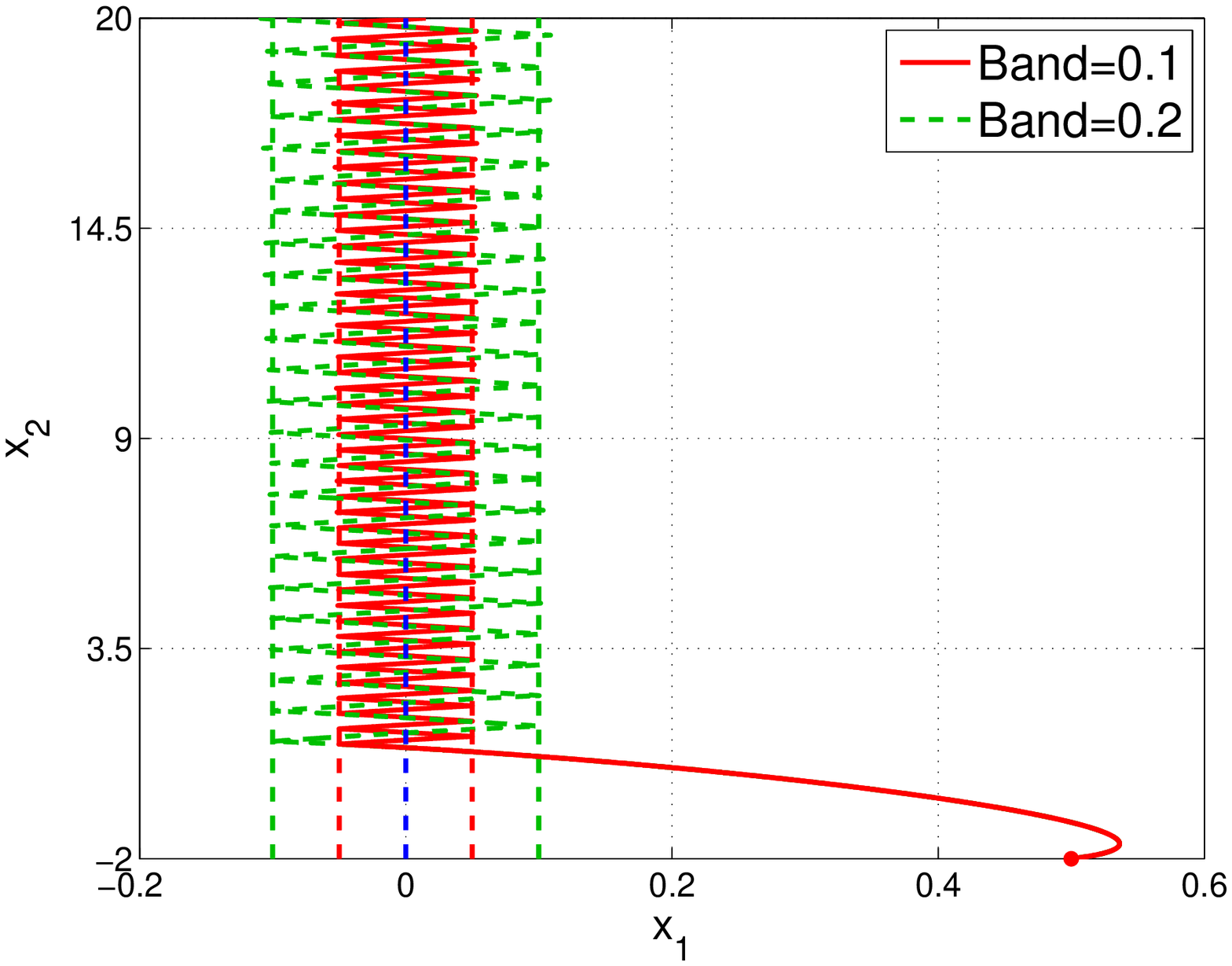} \\
		(a) & (b) \\
		\includegraphics[width=0.35\textwidth]{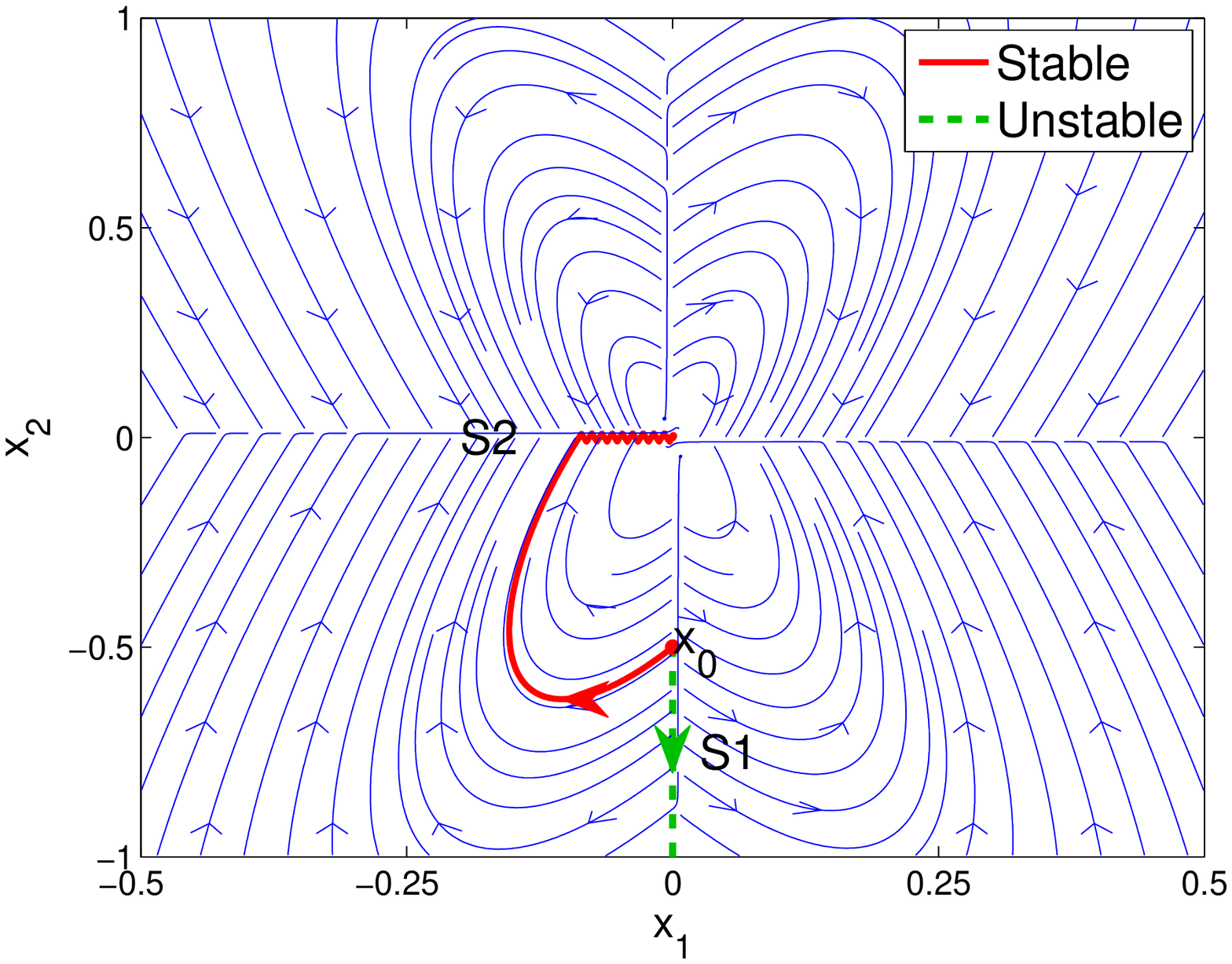} & 
		\includegraphics[width=0.35\textwidth]{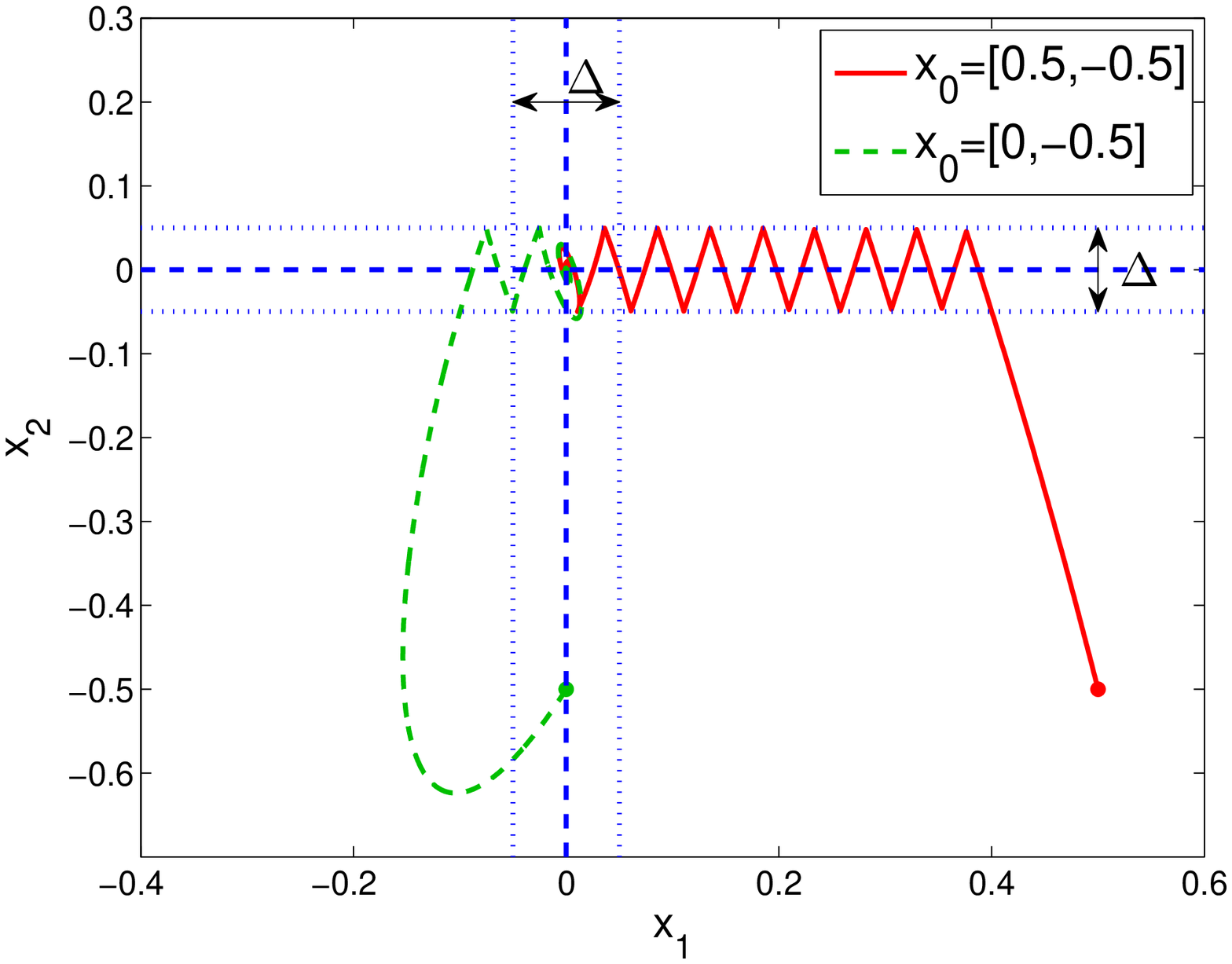} \\
		(c) & (d)
	\end{tabular}
	\caption{(a): The phase portrait of system (\ref{eq:unstablesm}), (b): The closed-loop trajectories of system (\ref{eq:unstablesm}) starting from $x_0=[0.5,-2]^T$, (c): The closed-loop trajectories of system (\ref{eq:nonuniquesm}) starting from $x_0=[0,-0.5]^T$, (d): The closed-loop trajectories of system (\ref{eq:nonuniquesm}) with hysteresis band $\Delta=0.1$.}
	\label{fig:unstablesm}
\end{figure*}

\begin{exmp} \textit{(Unstable Sliding Motion I)} Consider the following piecewise linear system:
\begin{align}
\dot{x}(t)=\left\{ \begin{matrix}
\begin{bmatrix}
-3 & -1 \\ 12 & 2
\end{bmatrix} x(t), & \textnormal{if } x_1(t) \ge 0 \\
\begin{bmatrix}
-3 & 1 \\ -12 & 2
\end{bmatrix} x(t), & \textnormal{if } x_1(t) \le 0 \end{matrix} \right.
\label{eq:unstablesm}
\end{align}
The closed-loop vector field is discontinuous on the surface $\{x \in \mathbb{R}^2: x_1=0\}$, resulting in a sliding motion. It is easy to see that both subsystems are stable and the closed-loop system will also be stable if no sliding motion occurs. However, the phase portrait in Fig.~\ref{fig:unstablesm}(a) indicates that the upper part of the switching surface (i.e. $\{x\in\R^2:x_1=0,x_2\ge 0\}$) is attractive\footnote{Roughly speaking, a sliding surface is attractive if all the trajectories starting from a sufficiently close neighborhood will converge to the surface. When the surface is differentiable, simple conditions are available to check its attractiveness~\cite{utkin1992sliding}.}and all the Filippov solutions will involve sliding motions regardless of the initial state locations. Fig.~\ref{fig:unstablesm}(b) shows two trajectories starting from $x_0=[0.5,-2]^T$ simulated with two hysteresis bands of different sizes. Both of them involve an unstable sliding motion that grows unbounded along the $x_2$ coordinate. 
\label{ex:unstablesliding}
\end{exmp}

Example \ref{ex:unstablesliding} indicates that excluding sliding motions may lead to an incorrect stability conclusion for switched systems. Therefore, sliding motions should be carefully considered during the switching law design. On the other hand, simply requiring stability for all sliding motions can be overly restrictive for many applications as illustrated in the next example.

\begin{exmp} \textit{(Unstable Sliding Motion II)} Consider the following piecewise linear system:
\begin{align}
\dot{x}(t)=\left\{ \begin{matrix}
\begin{bmatrix}
-3 & -1 \\ 12 & 2
\end{bmatrix} x(t), &  \begin{aligned} & \textnormal{if } x_1(t) \le 0, x_2(t) \ge 0 \\ & \textnormal{or } x_1(t) \ge 0, x_2(t) \le 0 \end{aligned} \\
\begin{bmatrix}
-3 & 1 \\ -12 & 2
\end{bmatrix} x(t), & \begin{aligned} & \textnormal{if } x_1(t) \ge 0, x_2(t) \ge 0 \\ & \textnormal{or } x_1(t) \le 0, x_2(t) \le 0 \end{aligned} \end{matrix} \right.
\label{eq:nonuniquesm}
\end{align}
As illustrated in Fig. \ref{fig:unstablesm}(c), the closed-loop vector field is discontinuous on the surface $S_1 \cup S_2$ where $S_1=\{x\in\R^2:x_1=0\}$ and $S_2=\{x\in\R^2:x_2=0\}$. According to the phase portrait in Fig.~\ref{fig:unstablesm}(c), if the initial state $x_0\in \R^2 \backslash S_1$, the trajectory will eventually reach $S_2$ followed by a stable sliding motion converging to the origin along $S_2$. If the initial state $x_0\in S_1$, the Filippov solutions are nonunique. Theoretically speaking, both stable and unstable Filippov solutions may exist (see Fig. \ref{fig:unstablesm}(c) for illustration). However, the sliding surface $S_1$, corresponding to the unstable sliding motion, is not attractive. As a result, any numerical simulation with nontrivial hysteresis bands can only produce stable trajectories as shown in Fig. \ref{fig:unstablesm}(d). In other words, the ideal unstable sliding motion (green dashed trajectory in Fig. \ref{fig:unstablesm}(c)) cannot occur in real systems. 
\label{ex:nonuniquesliding}
\end{exmp}

In Example \ref{ex:nonuniquesliding}, all the Filippov solutions converge to the origin except the sliding motions on $S_1$. Since the sliding surface $S_1$ is not attractive, sliding motions on this surface cannot sustain in reality and should not be considered in the analysis. Motivated by the above two examples, we will define two sets of Filippov solutions with one excluding sliding motions and the other including the sliding motions on attractive sliding surfaces (Example \ref{ex:unstablesliding}) but not the ones on nonattractive sliding surfaces (Example \ref{ex:nonuniquesliding}).

\begin{defn} Given an admissible switching law $\nu$ and a positive constant $\epsilon>0$, the set of Filippov solutions \textbf{excluding} sliding motions is defined as   
\begin{multline*}
\X_0(\epsilon,\nu) \triangleq \{x(\cdot;z,\nu): \|z\| \le \epsilon, \\ \dot{x}(t;z,\nu) = \fc(x(t;z,\nu)), \text{for } t \in \R_+ \text{ a.e.}\}
\end{multline*}
\end{defn}

Before introducing the set of Filippov solutions including sliding motions, we first establish the relation between the switching partitions and the vector fields of the active subsystems involved in the sliding motion by the following regularity condition. 

\begin{defn} Assume $x\in\partial\D$, $x$ is said to be regular if $\forall i \in \I_{sm}^a(x)$, there exists $\epsilon_i>0$ such that $x-\delta f_i(x)\in\D_i$, $\forall \delta \in (0,\epsilon_i)$.
\label{def:regular}
\end{defn}

Roughly speaking, if $x\in\partial\D_i$ is regular, then for any subsystem vector field $f_i(x)$ that contributes nontrivially to the velocity $\dot{x}$ (i.e. $\alpha_i(x)\neq 0$), we can conclude that $f_i(x)$ must point outwards $\D_i$. Note that if the switching boundaries $\partial \D$ are continuously differentiable, regularity condition can also be stated as $\forall i \in \I_{sm}^a(x)$, $\partial \phi_i(x)^T \cdot f_i(x)>0$ where $\phi_i(x)$ is a differentiable function such that $\D_i=\{x\in\R^n: \phi_i(x)<0\}$. Now we are ready to introduce the set of Filippov solutions that contains not only all the solutions in $\X_0$ but also the type of sliding motions in Example \ref{ex:unstablesliding}.

\begin{defn} Given an admissible switching law $\nu$ and a positive constant $\epsilon>0$, the set of Filippov solutions \textbf{including} sliding motions is defined as  
\begin{multline*} 
\X_{sm}(\epsilon,\nu) \triangleq \{x(\cdot;z,\nu): \|z\| \le \epsilon, \text{if } \exists (t_1,t_2), t_1<t_2 \\ \text{s.t. } \dot{x}(t;z,\nu) \neq \fc(x(t;z,\nu)) \text{ for } t\in (t_1,t_2) \text{ a.e.}, \\ \text{then } x(t;z,\nu) \text{ is regular } \forall t\in (t_1,t_2)\}
\end{multline*}
\end{defn}

We proceed to define asymptotic stability of the closed-loop system excluding and including sliding motions, respectively.

\begin{defn} System (\ref{eq:cl}) under an admissible switching law $\nu$ is asymptotically stable excluding sliding motions if (i) for any $\epsilon>0$, there exists $\delta>0$ such that for any Filippov solution $x(\cdot;z,\nu) \in \X_0(\delta,\nu)$, $\|x(t;z,\nu)\|<\epsilon, \forall t \in \R_+$. (ii) $\lim_{t\to\infty}x(t;z,\nu)=0$.
\end{defn}

\begin{defn} System (\ref{eq:cl}) under an admissible switching law $\nu$ is asymptotically stable including sliding motions if (i) for any $\epsilon>0$, there exists $\delta>0$ such that for any Filippov solution $x(\cdot;z,\nu) \in \X_{sm}(\delta,\nu)$, $\|x(t;z,\nu)\|<\epsilon, \forall t \in \R_+$. (ii) $\lim_{t\to\infty}x(t;z,\nu)=0$.
\label{def:stabsm}
\end{defn}

With these definitions, we are ready to formally state the switching stabilization problem to be studied in this paper.
\begin{prob}(Switching Stabilization Problem)
Find an admissible switching law $\nu$ under which the closed-loop system (\ref{eq:cl}) is asymptotically stable including sliding motions. 
\end{prob}
While there are many ways to study the switching stabilization problem, we are particularly interested in methods that (i) can lead to constructive ways to find a stabilizing switching law, (ii) allow for systematic stability analysis for sliding motions, and at the same time (iii) are less conservative and more general than existing methods. In the rest of this paper, we will first develop a control-Lyapunov function framework to accomplish the aforementioned objectives and then discuss its connections to the literature and its applications to switched nonlinear systems.

\section{A Piecewise Smooth Control-Lyapunov Function Framework}
\label{sec:framework}
The goal of this section is to develop a unified control-Lyapunov function (CLF) framework to solve the stabilization problems for switched nonlinear systems. We focus on piecewise smooth control-Lyapunov functions (PSCLFs), which allows for constructive design of stabilizing switching laws. To present the framework, we will first introduce the concept of PSCLFs, and then develop a strategy to construct switching laws using PSCLFs. Sufficient conditions will be derived to ensure closed-loop asymptotic stability including sliding motions. We will also discuss two important special cases of our PSCLF framework, namely, smooth CLFs and pointwise minimum CLFs, for which we will show that stable sliding motions are always guaranteed without any additional requirement on the sliding surfaces. 

\subsection{Piecewise Smooth Control-Lyapunov Function}
Nonsmooth CLFs have been well studied for stabilization of traditional nonlinear systems~\cite{S83,SS95,CLSS97,sontag1999stability}. Their relation to our switching stabilization problems has been discussed in the introduction section. This subsection will introduce a particular class of nonsmooth CLFs, namely, piecewise smooth CLFs, for switching stabilization problems. We first give a formal definition for piecewise smooth functions.
\begin{defn} \label{def:pwsf} \textit{(Piecewise Smooth Function)} A function $g:\mathbb{R}^n \to \mathbb{R}$ is called piecewise smooth if it is continuous and there exists a finite collection of disjoint and open sets, $\Omega_1, \ldots, \Omega_m\subseteq \R^n$ such that (i) $\cup_{j\in\N_m}\bar\Omega_j = \R^n$, (ii) $g$ is continuously differentiable on $\Omega_j$, $\forall j\in\N_m$, and (iii) $\partial\Omega_j$ is a differential manifold for each $j\in\N_m$.
\end{defn}

\begin{rem}
Piecewise smooth functions can be defined in different ways~\cite{LS93, subbotin1995generalized}. In our definition, each partition $\Omega_j$ is not required to be connected. In addition, since each $\Omega_j$ is a full-dimensional open set in $\R^n$, its boundary $\partial\Omega_j$ must have measure zero, i.e. $\mu(\partial \Omega_j)=0$. 
\end{rem}

An important property of piecewise smooth function is the existence of directional derivative anywhere in the state space as stated in the following lemma.

\begin{lem} \label{lem:dd} (Directional Derivative \cite{subbotin1995generalized}) A piecewise smooth function $g:\mathbb{R}^n \to \mathbb{R}$ as defined in Definition~\ref{def:pwsf} is directionally differentiable on $\mathbb{R}^n$, i.e., the limit $Dg(x;\eta) \triangleq \lim_{\delta \downarrow 0} \frac{1}{\delta}(g(x+\delta \eta)-g(x))$ exists, $\forall x, \eta \in \mathbb{R}^n$. 
\end{lem}

Control-Lyapunov function is a useful tool to study stabilization problems. To be less conservative, we will focus on piecewise smooth control-Lyapunov functions for general switched nonlinear systems (\ref{eq:ol}) defined as follows.

\begin{defn} \label{def:clf} (Piecewise Smooth Control-Lyapunov Function (PSCLF)) A piecewise smooth function $V: \mathbb{R}^n \to \mathbb{R}_+$ is called a PSCLF if there exists another continuous function $W: \mathbb{R}^n \to \mathbb{R}_+$ such that the following conditions hold:
\begin{align}
& V(x)>0, W(x)>0, \forall x \neq 0, \quad V(0)=0 \label{cond:pd} \\
& \mathcal{L}_{\beta}=\{x:V(x) \le \beta\} \text{ is bounded for each } \beta \label{cond:prop} \\
& \min_{i \in \mathcal{Q}} DV(x;f_i(x)) \le -W(x), \quad \forall x \in \mathbb{R}^n  \label{cond:inf} 
\end{align}
In addition, the pair of functions $(V,W)$ is called a CLF pair.
\end{defn}

\begin{rem} The PSCLF defined above can be viewed as a special class of nonsmooth CLFs~\cite{SS95}. Roughly speaking, the existence of such a function implies the stabilizability of the system in the sample-and-hold sense~\cite{CLSS97}. The main complication for switching stabilization lies in how to constructively find a stabilizing switching law in the Filippov sense and how to guarantee the closed-loop stability including the sliding motions on the attractive sliding surfaces.
\end{rem}

We refer to (\ref{cond:pd}) as the positive definite condition, refer to (\ref{cond:prop}) as the radially unbounded condition, and refer to (\ref{cond:inf}) as the decreasing condition. Note that the decreasing condition is given in terms of directional derivative, which is well-defined for piecewise smooth functions according to Lemma~\ref{lem:dd}. In the rest of this section, we will first develop a constructive way to design a switching law using a given PSCLF, and then carefully study various stability properties of the closed-loop system. 

\subsection{Switching Law Construction}
If system (\ref{eq:ol}) admits a PSCLF $V$, then it can be used to construct a stabilizing switching law. The main idea is to select the subsystem along whose vector field the PSCLF decreases at the fastest rate. Naturally, one may want to directly construct the switching law as 
\begin{align}
\hat\nu(x) = \argmin_{i\in\Q} DV(x;f_i(x)), \quad \forall x\in\R^n. \label{eq:nuhat}
\end{align}
The switching law $\hat\nu:\R^n\to 2^{\Q}$ defined above is set-valued and the part of the state space corresponding to multiple minimizers (i.e. $X_s\triangleq \{x\in\R^n: |\hat\nu(x)|>1\}$) is where sliding motions may possibly occur. To understand the stability behaviors of all the closed-loop Filippov solutions, there are several technical issues requiring special attention. First, the set $X_s$ may not have measure zero, which complicates the derivation and analysis of the Filippov solutions of the closed-loop system. In addition, even if $X_s$ is of measure zero, it may intersect the nonsmooth surface of $V$, introducing additional challenges in analyzing the sliding motions. Some aspects of these technical issues regarding the switching law $\hat{\nu}$ are further discussed in Example \ref{ex:DV} in Section 3.3. To better address these issues, we will develop a slightly different method to construct the switching law, which can facilitate our discussion on Filippov solutions and their stabilities. 

We first introduce some notations. Suppose we are given a general PSCLF $V:\R^n\to\R_+$ with partition sets $\{\Omega_j\}_{j=1}^m$. According to Definition~\ref{def:pwsf}, each partition $\Omega_j$ is an open set in $\R^n$ with boundary $\partial\Omega_j$ and closure $\bar\Omega_j$, $j\in\N_m$. The union of all the partition boundaries will be denoted by $\partial\Omega\triangleq \cup_{j\in\N_m} \partial\Omega_j$. Let $J(x)$ be the set of partitions whose closure $x$ belongs to, i.e., 
\begin{align}
J(x) = \{j\in\N_m: x\in \bar\Omega_j\},  \quad x\in\R^n. \label{eq:J}
\end{align}
$J(x)$ is a singleton if $x$ lies inside one of the partitions, and is set-valued if $x\in\partial\Omega$. Denote by $V_j:\bar\Omega_j\to\R$ the restriction of the function $V$ to the closure $\bar\Omega_j\subset\R^n$. By Definition~\ref{def:pwsf}, $V_j$ is continuously differentiable on $\Omega_j$. Therefore, we have
\begin{align}
DV_j(x;\eta) = <\nabla V_j(x), \eta>, \forall x\in \Omega_j \text{ and } \eta\in \R^n, \label{eq:DVjsm}
\end{align}
where $\nabla V(x)$ denotes the gradient of $V$ at $x$. For any boundary point $\hat x\in \partial\Omega_j$, we define $\nabla V_j(\hat x) \triangleq \lim_{x\in\Omega_j, x\to\hat x} \nabla V_j(x)$. Since $V_j$ is continuously differentiable on $\Omega_j$, we have
\begin{multline}
V_j(z) - V_j(x) = <\nabla V_j(x), z-x> + o\left(\|z-x\|\right), \\ \forall x, z\in \bar\Omega_j, j\in\N_m, \text{where }\lim_{z\to x}\frac{o\left(\|z-x\|\right)}{\|z-x\|}=0.
\label{eq:dvj}
\end{multline}
With these notations, the switching law can be constructed as follows:
\begin{enumerate}
\item For each $j\in\N_m, i\in\Q$, define 
\begin{multline}
\D_{j,i}\triangleq \Big\{x\in\Omega_j:\langle \nabla V_j(x), f_i(x) \rangle < \\ \min_{q \in \Q,q \neq i} \langle \nabla V_j(x), f_q(x)\rangle \Big\};
\label{eq:Dij}
\end{multline}
\item For each $i\in\Q$, define $\D_i \triangleq \mathop{\cup}_{j\in\N_m} \D_{j,i}$ and $\partial\D \triangleq \cup_{i\in\Q} \partial\D_i$;
\item Construct the switching law as: 
\begin{align}
\hspace{-1cm}
\nu(x) = \left\{ \begin{array}{l}
i, \quad \text{if }\exists i\in\Q, \text{ s.t. } x\in \D_i \\
\min\{i\in\Q: x\in\partial \D_i\}, \quad \text{if } x\in \partial\D
\end{array} \right.
\label{eq:nu}
\end{align}
\end{enumerate}
The above procedure will be represented by an operator, denoted by $\S$, which maps a PSCLF $V$ to a particular switching law $\nu = \S[V]$. Due to the continuities of all the vector fields and $\nabla V_j$, the set $\D_{j,i}$ defined in~(\ref{eq:Dij}) is an open set in $\R^n$ for each $j\in\N_m$ and $i\in\Q$. Hence, the set $\D_i$ is also open for each $i\in\Q$. In addition, the boundaries $\partial \D_i$, $i\in\Q$, and their union $\partial \D$ are all of measure zero. In summary, the sets $\{\D_i\}_{i\in\Q}$ are open and disjoint. Therefore, $\nu$ defined in (\ref{eq:nu}) is an admissible switching law with switching partitions $\{\D_i\}_{i\in\Q}$ and switching boundaries $\{\partial\D_i\}_{i\in\Q}$.

It is also easy to see that the switching law $\nu$ defined in~(\ref{eq:nu}) assigns each state $x\in\R^n$ to a unique subsystem. It is worth mentioning that there are many different ways to design the switching laws on the boundary set $\partial\D$. All of them will lead to the same closed-loop Filippov set-valued map as $\partial\D$ is of measure zero. The particular construction given in~(\ref{eq:nu}) ensures that $\nu(\hat{x})$ at any boundary point $\hat{x}\in\partial \D$ is a continuous extension of $\nu$ on one of the partitions whose boundaries touch $\hat{x}$. 

Next, we evaluate the closed-loop system (\ref{eq:cl}) under the switching law $\nu$ defined in (\ref{eq:nu}). Clearly, the closed-loop vector field $\fc$ is measurable and locally essentially bounded. Thus, Filippov solutions always exist. In addition, $\fc$ is continuous inside each switching partition $\D_i$ and piecewise continuous around any boundary point $\hat{x}\in\partial\D$, thus the Filippov set-valued map is given by \cite{C09}:
\begin{align}
F\left[\fc\right](x) = \left\{ \begin{array}{l}
\{f_i(x)\}, \quad \text{if }\exists i\in\Q, \text{ s.t. } x\in\D_i \\
\co{f_i(x): i\in\Q, x\in\bar\D_i}, \text{ o/w}
\end{array} \right.
\label{eq:Fv}
\end{align}
The rest of this section will analyze and discuss various stability properties of the closed-loop Filippov solutions with Filippov set-valued map $F\left[ \fc \right]$ defined in~(\ref{eq:Fv}).

\subsection{Piecewise Smooth Control-Lyapunov Function Theorems}
In this subsection, we will establish conditions to ensure stability for the closed-loop system. The main result consists of two parts. First, we will show that if $V$ is a PSCLF, then the closed-loop system~(\ref{eq:cl}) under the switching law $\nu = \S[V]$ is asymptotically stable provided there is no sliding motion. Second, we will derive an additional condition for the PSCLF $V$ to guarantee closed-loop stability including sliding motions as defined in Definition~\ref{def:stabsm}.

Throughout this subsection, we will assume $(V, W)$ is a CLF pair, where $W:\R^n\to\R_+$ is a continuous nonnegative function, and $V:\R^n\to\R_+$ is a PSCLF with partitions $\{\Omega_j\}_{j\in\N_m}$ and nonsmooth boundaries $\partial\Omega = \cup_{j\in\N_m} \partial\Omega_j$. In addition, the switching law $\nu =\S[V]$ is generated from $V$ as described in the previous subsection with switching partitions $\{\D_i\}_{i\in\Q}$ and switching boundaries $\partial \D=\cup_{i\in\Q} \partial\D_i$. 

Roughly speaking, the key of our stability analysis is to ensure the PSCLF $V$ decreases along any closed-loop trajectory $x(t)$. If there is no sliding motion, we have $\dot x(t) = f_{\nu(x(t))}(x(t))$, for almost all $t\in\R_+$. In this case, we need to check the directional derivative along the closed-loop vector field, i.e. $DV(x;f_{\nu(x)}(x))$, $x\in\R^n$. On the other hand, if there is a sliding motion, then there exists a nontrivial time interval $(t_1,t_2)$ such that $\dot x(t) =\sum_{i\in\I_{sm}(x(t))}\alpha_i(x(t))f_i(x(t))$, for almost all $t\in(t_1,t_2)$. In this case, merely looking at $DV(x,f_{\nu(x)}(x))$ is no longer enough; we need to guarantee that $V$ decreases along the switching boundaries $\partial \D$. In other words, we need to further check $DV(x;\sum_{i\in\I_{sm}(x)}\alpha_i(x)f_i(x))$ for $x \in \partial\D$. Before presenting the main stability results, we first derive some useful properties for these two types of directional derivatives. 

\begin{lem}
\label{lem:DV}
For any point not on the nonsmooth boundaries, i.e. $x\notin\partial\Omega$, we have 
\begin{align*}
DV(x;f_{\nu(x)}(x)) =\min_{i \in \Q} DV(x;f_i(x))\le -W(x).
\end{align*}
\end{lem}
\begin{pf}
Assume $x \in \Omega_j$ for some $j \in \N_m$. According to the property of $V$ in (\ref{eq:DVjsm}) and the construction of $\nu$ in (\ref{eq:Dij}), $DV(x;f_{\nu(x)}(x)) = \langle \nabla V_j(x), f_{\nu(x)}(x) \rangle \le \langle \nabla V_j(x), f_i(x) \rangle=DV_j(x;f_i(x))$, $\forall i \in \Q$. In other words, $DV(x;f_{\nu(x)}(x)) = \min_{i \in \mathcal{Q}}DV_j(x;f_i(x))$. It follows from the decreasing condition of $V$ in (\ref{cond:inf}) that $\min_{i \in \mathcal{Q}}DV_j(x;f_i(x)) \le -W(x)$, $\forall x \in \Omega_j$. Therefore, we conclude that $DV(x;f_{\nu(x)}(x)) \le -W(x)$, $\forall x \in \Omega_j$ for some $j \in \N_m$. \qed
\end{pf}

As mentioned at the beginning of Section 3.2, one of the reasons to use the switching law $\nu$ defined in (\ref{eq:nu}) instead of the more natural choice $\hat{\nu}$ defined in (\ref{eq:nuhat}) is that the latter one introduces additional challenges in analyzing the sliding motions on the nonsmooth surface. In general, the two switching laws give different results on the nonsmooth surface. In particular, on the nonsmooth surface, the directional derivative along closed-loop trajectory under $\nu$ differs from that under $\hat{\nu}$, i.e., $DV(x;f_{\nu(x)}(x))\neq DV(x;f_{\hat{\nu}(x)}(x))$, and the set of subsystems that achieve the minimal decreasing rate under $\nu$ differs from the one under $\hat{\nu}$, i.e., $\I_{sm}(x)\neq \hat{\nu}(x)$. To better illustrate the subtle differences between the two switching laws on the nonsmooth surface, we give the following example.

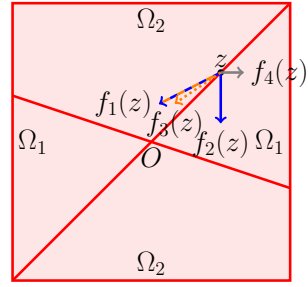
\begin{figure}[ht!]
\centering
\resizebox{0.5\linewidth}{!}{ \begin{tikzpicture}
\filldraw [ultra thick,red,fill opacity=0.1] (0,0) -- (6,0) -- (6,6) -- (0,6) -- cycle;

\draw [ultra thick,red] (0,4) -- (6,2);
\draw [ultra thick,red] (0,0) -- (6,6);

\coordinate (p1) at (3,6);
\coordinate (p2) at (3,0);
\coordinate (p3) at (0,3);
\coordinate (p4) at (6,3);

\node [fill=none,draw=none,below] at (p1) {\Large $\Omega_2$};
\node [fill=none,draw=none,above] at (p2) {\Large $\Omega_2$};
\node [fill=none,draw=none,right] at (p3) {\Large $\Omega_1$};
\node [fill=none,draw=none,left] at (p4) {\Large $\Omega_1$};

\coordinate (state) at (4.5,4.5);
\fill (state) circle (2pt) node[above] {\Large $z$};
\node [fill=none,draw=none,below] at (3,3) {\Large $O$}; 


\draw [ultra thick,->,blue,dash pattern= on 5pt off 5pt] (state) -- +(207:1.5cm) node[black,left] {\Large $f_1(z)$};
\draw [ultra thick,orange,dash pattern= on 5pt off 5pt,dash phase=5pt] (state) -- +(207:1.5cm);
\draw [ultra thick,->,blue] (state) -- +(-90:1.1cm) node[black,below] {\Large $f_2(z)$};
\draw [ultra thick,->,orange,dotted] (state) -- +(214:1.2cm) node[black,below] {\Large $f_3(z)$};
\draw [ultra thick,->,gray] (state) -- +(0:0.5cm) node[black,right] {\Large $f_4(z)$};

\end{tikzpicture} }
\caption{Subsystem vector fields on the nonsmooth surface of Example~\ref{ex:DV}}
\label{fig:exDV}
\end{figure}

\begin{exmp}
Consider a switched linear system $\dot{x}=A_ix$, $i \in \mathcal{Q}=\{1, \cdots, 4\}$ and given a PSCLF $V=x^TP_jx, \forall x \in \Omega_j$, $j \in \N_2$, where 
\begin{align*}
& A_1=\begin{bmatrix} 0 & -4 \\ -2 & 0 \end{bmatrix},
A_2=\begin{bmatrix} 0 & 0 \\ -\frac{10}{3} & 0 \end{bmatrix},
A_3=\begin{bmatrix} 0 & -3 \\ -2 & 0 \end{bmatrix}, \\
& A_4=\begin{bmatrix} 0 & 1 \\ 1 & -1 \end{bmatrix},
P_1=\begin{bmatrix} 2 & 0 \\ 0 & 1 \end{bmatrix},
P_2=\begin{bmatrix} 1 & -1 \\ -1 & 4 \end{bmatrix}.
\end{align*}
The nonsmooth surface of $V$ is given by $\{x \in \mathbb{R}^2: x^T(P_1-P_2)x=0\}=\{x\in\R^2:x_1-x_2=0\} \cup \{x\in\R^2:x_1+3x_2=0\}$, which consists of two lines that intersect at the origin. In this example, $DV_j(x;f_i(x))$, $j \in \N_2$, $i \in \mathcal{Q}$ takes the form of $DV_j(x;f_i(x))=x^T(A_i^TP_j+P_jA_i)x$. We pick a point of the form $z=[c,c]^T \in \partial \Omega$ and analyze its directional derivatives. Partitions $\{\Omega_j\}_{j\in\N_2}$ and subsystem vector fields at $z$ are shown in Fig.~\ref{fig:exDV}. The directional derivatives at $z$ can be easily calculated as
\begin{align}
\begin{aligned}
& DV_1(z;f_1(z))=-20c^2,DV_1(z;f_2(z))=-\frac{20}{3}c^2, \\
& DV_1(z;f_3(z))=-16c^2,DV_1(z;f_4(z))=4c^2, \\
& DV_2(z;f_1(z))=-12c^2,DV_2(z;f_2(z))=-20c^2, \\
& DV_2(z;f_3(z))=-12c^2,DV_2(z;f_4(z))=0.
\end{aligned}
\label{eq:exdv}
\end{align}
The above results indicate the following relations among the directional derivatives:
\begin{align}
\begin{aligned}
& DV_1(z;f_1(z))=DV_2(z;f_2(z)) \\
& DV_1(z;f_1(z)) < DV_1(z;f_i(z)), \quad \forall i \neq 1 \\
& DV_2(z;f_2(z)) < DV_2(z;f_i(z)), \quad \forall i \neq 2
\end{aligned}
\label{eq:exdvsmy}
\end{align}
According to the switching law $\nu$ defined in~(\ref{eq:nu}), the minimal directional derivative at $z$ is 
\begin{align*}
\min_{j\in\N_2}\min_{i\in\Q}DV_j(z;f_i(z))=-20c^2,
\end{align*}
which is actually the minimum of the eight terms in (\ref{eq:exdv}), and the corresponding set of subsystems that achieve this minimum is $\I_{sm}(z)=\{1,2\}$ (solid and dashed blue arrows in Fig.~\ref{fig:exDV}). However, if using the commonly used switching law $\hat{\nu}$ defined in (\ref{eq:nuhat}), the minimal directional derivative at $z$ is
\begin{multline*}
\min_{i\in\Q}DV(z;f_i(z))=\min\{DV_1(z;f_2(z)), \\ DV_1(z;f_4(z)),DV_2(z;f_1(z)),DV_2(z;f_3(z))\}=-12c^2
\end{multline*}
where the first equality follows from the definition of directional derivative. The set of subsystems that achieve this minimum is $\hat{\nu}(z)=\{1,3\}$ (dotted and dashed orange arrows in Fig.~\ref{fig:exDV}). In summary, the switching law $\nu$ defined in (\ref{eq:nu}) gives a different result from the switching law $\hat{\nu}$ on the nonsmooth surface where PSCLF $V$ decreases faster under $\nu$ than under $\hat{\nu}$.
\label{ex:DV}
\end{exmp}

In general, $DV(x;f_{\nu(x)}(x))\neq\min_iDV(x;f_i(x))$ on $\partial\Omega$ as illustrated in Example \ref{ex:DV}. The expression for $DV(x;f_{\nu(x)}(x))$ on the nonsmooth boundary $\partial\Omega$ can be quite involved. To address this issue, we introduce a set-valued map defined below.

\begin{defn} \label{def:Mj} For any $j\in\N_m$, the set-valued map $\mathcal{M}_j:\partial \Omega_j \to 2^{\Q}$ is defined as
\begin{multline}
\mathcal{M}_j(x)=\{i \in \mathcal{Q}: \exists \{z_k\}_{k\ge 0}\in \Omega_j\cap \D_i \\ \text{such that } \lim_{k \to \infty} z_k = x \}, x\in\partial\Omega_j.
\label{eq:M}
\end{multline}
\end{defn}

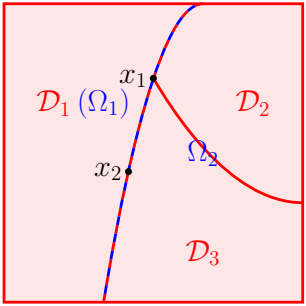
\begin{figure}[ht!]
\centering
\resizebox{0.5\linewidth}{!}{ \begin{tikzpicture}
\filldraw [ultra thick,red,fill opacity=0.1] (0,0) -- (6,0) -- (6,6) -- (0,6) -- cycle;

\draw [ultra thick,red,dash pattern= on 5pt off 5pt] (4,6) parabola (2,0);
\draw [ultra thick,blue,dash pattern= on 5pt off 5pt,dash phase=5pt] (4,6) parabola (2,0);
\draw [ultra thick,red] (6,2) parabola (3,4.5);

\coordinate (x1) at (3,4.5);
\coordinate (x2) at (2.5,2.625);
\fill (x1) circle (2pt) node[left] {\Large $x_1$};
\fill (x2) circle (2pt) node[left] {\Large $x_2$};

\node [fill=none,draw=none,color=red] at (1,4) {\Large $\mathcal{D}_1$}; 
\node [fill=none,draw=none,color=red] at (5,4) {\Large $\mathcal{D}_2$}; 
\node [fill=none,draw=none,color=red] at (4,1) {\Large $\mathcal{D}_3$}; 

\node [fill=none,draw=none,color=blue] at (2,4) {\Large $(\Omega_1)$}; 
\node [fill=none,draw=none,color=blue] at (4,3) {\Large $\Omega_2$}; 

\end{tikzpicture} }
\caption{Illustrating example for $\M_j(x)$}
\label{fig:exM}
\end{figure}

Note that in the above definition $\D_i$ and $\Omega_j$ are both open sets, so is their intersection. The above definition indicates that $x\in\partial\Omega_j$ is a limit point of $\Omega_j\cap\D_i$ for all $i\in\M_j(x)$. Roughly speaking, $\M_j(x)$ contains the set of subsystem indices that are chosen by the switching law $\nu$ within an arbitrarily small neighborhood $\mathcal{N}(x;\epsilon)\cap\Omega_j$, for all sufficiently small $\epsilon>0$. For example, consider the partitions shown in Fig.~\ref{fig:exM}, where $\Omega_1$ is the same as $\D_1$ and $\Omega_2$ is the union of $\D_2$ and $\D_3$. For this particular example, we have $\M_1(x_1)=\{1\}$, $\M_2(x_1)=\{2,3\}$, $\M_1(x_2)=\{1\}$, $\M_2(x_2)=\{3\}$.

The following properties of $\M_j$ are helpful for deriving the main stability results.

\begin{lem}
\label{lem:Mj}
Let $x\in \partial\Omega$ and let $J(x)$ be defined in~(\ref{eq:J}). 
\begin{enumerate}[label=\roman*)]
\item If $x\in\D_i$ for some $i\in\Q$, then $\M_j(x) = \{i\}$, $\forall j\in J(x)$. 
\item If $x\in\partial\D$, then $\cup_{j\in J(x)} \M_j(x) = \I_{sm}(x)$. 
\end{enumerate}
\end{lem}

\begin{pf}
Part i) of this lemma follows immediately from Definition~\ref{def:Mj} and the fact that the sets $\{\D_i\}_{i\in\Q}$ are open and mutually disjoint. For part ii), it suffices to show that if $i\in\I_{sm}(x)$, then $i\in\M_j(x)$ for some $j\in J(x)$. We thus fix an arbitrary subsystem index $i\in \I_{sm}(x)$. Clearly, $x\in\partial\D_i$ and $\exists \delta_1>0$, such that $\mathcal{N}(x,\epsilon)\cap \D_i$ is a nonempty open set  $\forall \epsilon\in(0,\delta_1)$. In addition, according to the definition of $J(x)$, we know there must exist $\delta_2>0$ such that $\mathcal{N}(x,\epsilon) = \cup_{j\in J(x)} \left(\bar\Omega_j \cap \mathcal{N}(x,\epsilon)\right)$, $\forall \epsilon\in (0, \delta_2)$. Therefore, $\exists \hat j\in J(x)$ such that $\mathcal{N}(x,\epsilon)\cap \D_i\cap\bar \Omega_{\hat j} \neq \emptyset$, which in turn implies $\mathcal{N}(x,\epsilon)\cap \D_i\cap\Omega_{\hat j}$ is a nonempty open set for all $\epsilon\in (0,\delta_0)$, where $\delta_0=\min\{\delta_1,\delta_2\}$. Let $\epsilon_k = \min\{\frac{1}{2^k}, \delta_0\}$ and pick any $z_k\in \mathcal{N}(x, \epsilon_k)\cap \D_i\cap\Omega_{\hat j}$. Clearly, $z_k\to x$, as $k\to\infty$. Therefore, $i\in \M_{\hat j}(x)$, which completes the proof. \qed
\end{pf}

Some properties of directional derivatives on the nonsmooth surface can be revealed with the help of the set-valued map $\M_j$ defined in Definition \ref{def:Mj}.

\begin{lem}
\label{lem:DVj}
For any $j\in\N_m$, $x\in\partial\Omega_j$, and $q\in\M_j(x)$, we have
\begin{enumerate}[label=\roman*)]
\item $DV_j(x;f_q(x))=\min_{i \in \Q} DV_j(x;f_i(x))$;
\item $DV_j(x;f_q(x))\le -W(x)$.
\end{enumerate}
\end{lem}

\begin{pf}
We fix $x \in \partial \Omega_j$ for some $j \in \N_m$ and $q \in \M_j(x)$. By Definition \ref{def:Mj}, $\exists \{z_k\}_{k\ge 0}\in \Omega_j\cap \D_q$ such that $\lim_{k \to \infty}z_k=x$. Based on the definition of $\nu$ in (\ref{eq:nu}), $\nu(z_k)=q$, which implies $DV_j(z_k;f_{\nu(z_k)}(z_k))=DV_j(z_k;f_q(z_k))$. Apply Lemma \ref{lem:DV} to $\{z_k\}_{k \ge 0} \in \Omega_j$, $DV_j(z_k;f_{\nu(z_k)}(z_k))=\min_{i \in \Q} DV_j(z_k;f_i(z_k))\le -W(z_k)$. From above, we have the following relation.
\begin{align}
\hspace{-0.1cm}
DV_j(z_k;f_q(z_k))=\min_{i \in \Q} DV_j(z_k;f_i(z_k))\le -W(z_k) \label{eq:DVjzk}
\end{align}
Due to the continuity of $DV_j(\cdot,f_q(\cdot))$ for all $j \in \N_m$, $q \in \Q$ and the continuity of $\min_i g_i$ for finite number of continuous function $g_i$, the equality part in (\ref{eq:DVjzk}) still holds in the limit as $k \to \infty$, which completes the proof of part i). Part ii) directly follows from the inequality part in (\ref{eq:DVjzk}) and the continuity of $W(\cdot)$. \qed
\end{pf}

Although it is not generally true for $x\in\partial\Omega$ that $DV(x;f_{\nu(x)}(x))=\min_{i\in\Q}DV(x;f_i(x))$ (see Example \ref{ex:DV}), a similar conclusion can still be obtained if we look at the part of nonsmooth surface that does not contain switching boundaries.

\begin{lem}
\label{lem:DVboundary}
For any $x\in\partial\Omega \backslash \partial\D$, we have
\begin{align*}
DV(x,f_{\nu(x)}(x))\le -W(x).
\end{align*}
\end{lem}

\begin{pf}
Assume $x \in \D_q$ for some $q \in \Q$. According to part i) of Lemma \ref{lem:Mj}, $q \in \M_j(x)$, $\forall j \in J(x)$. By the definition of $J(x)$, $x\in \partial \Omega_j$, $\forall j\in J(x)$. By part ii) of Lemma \ref{lem:DVj}, $DV_j(x;f_q(x))\le -W(x)$, $\forall j \in J(x)$. By the definition of directional derivative, for any $\eta\in\R^n$ and $x\in\partial\Omega$, $DV(x;\eta)=DV_{j^*}(x;\eta)$ for some $j^*\in J(x)$. The desired result then follows from $DV(x;f_{\nu(x)}(x))=DV_{j^*}(x;f_{\nu(x)}(x))$ for some $j^* \in J(x)$. \qed
\end{pf}

With the above lemmas, we are now ready to state our first stability result.

\begin{thm} \label{thm:pwsclf1} (PSCLF Theorem I) The closed-loop system~(\ref{eq:cl}) under the switching law $\nu = \S[V]$ is globally asymptotically stable excluding sliding motions. 
\end{thm}

\begin{pf}
Let $x(t)$ be the closed-loop trajectory starting from an arbitrary initial state $x_0\in\R^n$ under the switching law $\nu=\S[V]$. Since there is no sliding motion, we have $x(t) \not\in \partial \D$ for almost all $t\in\R_+$. Thus, we have $\frac{d}{dt}V(x(t))=DV\left(x(t);f_{\nu(x(t))}(x(t))\right)\le -W(x(t))$, for almost all $t \in \R_+$, where the last inequality is due to Lemma~\ref{lem:DV} and Lemma~\ref{lem:DVboundary}. The rest of the proof follows directly from the classical Lyapunov theorem proof (\cite[Theorem 4.1]{khalil2002nonlinear}) by replacing the Lie derivative with the directional derivative $DV(x;\dot{x})$. \qed
\end{pf}

The next goal is to analyze the stability of sliding motions that may occur along a closed-loop trajectory. We first introduce another technical lemma regarding the directional derivative along a tangent direction of the nonsmooth boundary. 

\begin{lem} \label{lem:ddBoundary} For any $x\in\partial\Omega$ and $j, j^\prime\in\N_m$ satisfying $\mathcal{T}_x(\partial\Omega_j)\cap \mathcal{T}_x(\partial\Omega_{j^\prime})\neq \emptyset$, we have $DV(x;\eta)=DV_j(x;\eta) =DV_{j'}(x;\eta)$, $\forall\eta\in\mathcal{T}_x(\partial\Omega_j)\cap \mathcal{T}_x(\partial\Omega_{j^\prime})$.
\end{lem}

\begin{pf} We pick an arbitrary point $x \in \partial \Omega$ and let $j,j' \in \N_m$ satisfy $\mathcal{T}_x(\partial\Omega_j)\cap \mathcal{T}_x(\partial\Omega_{j^\prime})\neq \emptyset$. It is easy to see that $j,j' \in J(x)$. We then pick an arbitrary vector $\eta \in \mathcal{T}_x(\partial\Omega_j)\cap \mathcal{T}_x(\partial\Omega_{j^\prime})$. According to the definition of tangent vector, for any $\delta>0$, there exists a function $h:\R\to\R^n$ with $\lim_{\delta\to 0}\|h(\delta)\|=0$ such that $x+(\eta+h(\delta))\delta \in \bar{\Omega}_j$. By the definition of $J(x)$ in (\ref{eq:J}), we have $x\in\bar{\Omega}_j$. It follows from (\ref{eq:dvj}) that $V_j(x+(\eta+h(\delta))\delta)-V_j(x)=\langle \nabla V_j(x),(\eta+h(\delta))\delta \rangle+o(\|\eta+h(\delta)\|\delta)$ where $\lim_{\delta\to 0} \frac{o(\|\eta+h(\delta)\|)}{\|\eta+h(\delta)\|}=0$. By first dividing $\delta$ and then taking the limit of $\delta\downarrow 0$ on both sides of the equation above, we have
\begin{align}
\hspace{-0.2cm}
\lim_{\delta \downarrow 0} \frac{1}{\delta} (V(x+(\eta+h(\delta))\delta)-V(x))=\langle \nabla V_j(x),\eta \rangle. \label{eq:dd}
\end{align}
Since $V$ satisfies Lipschitz condition, the limit in (\ref{eq:dd}) coincides with $DV(x;\eta)$ and hence $DV(x;\eta)=DV_j(x;\eta)$. The above argument also holds for $j'$, which completes the proof. \qed
\end{pf} 

\begin{thm} \label{thm:pwsclf2} (PSCLF Theorem II) The closed-loop system~(\ref{eq:cl}) under the switching law $\nu=\S[V]$ is globally asymptotically stable including sliding motions if for any $x \in \partial \Omega \cap \partial \mathcal{D}$, there exists $j\in J(x)$ such that 
\begin{align}
DV_j(x;f_q(x)) \le -W(x), \forall q \in \M_j^c(x) \cap \I_{sm}^a(x).
\label{eq:bdcond}
\end{align}
\end{thm}

\begin{pf}
Let $x(t)$ be the closed-loop trajectory of system (\ref{eq:cl}) under the switching law $\nu=\S[V]$. We want to show that $DV(x(t);\dot{x}(t)) \le -W(x(t))$ for $t \in \mathbb{R}_+$ a.e.. By Theorem \ref{thm:pwsclf1}, we only need to show the inequality for the part of the trajectory that involves sliding motions. Let $(t_1,t_2)$ be the time interval that sliding motion occurs, i.e., $\dot{x}(t)=\sum_{i\in\I_{sm}(x(t))}\alpha_i(x(t))f_i(x(t))$ where $\sum_{i\in\I_{sm}(x(t))}\alpha_i(x(t))=1, t\in (t_1,t_2)$ a.e.. The proof is divided into two cases: i) the ones on the smooth region ($\{\Omega_j\}_{j\in\N_m}$) and ii) the ones on the nonsmooth boundary ($\partial\Omega$). We want to show that for both cases $V$ satisfies 
\begin{align}
DV(x(t);\dot{x}(t)) \le -W(x(t)) \text{ for } t \in (t_1,t_2) \text{ a.e..} \label{eq:psclfcond}
\end{align}
To simplify notation, we will use $x$ instead of $x(t)$ in the rest of the proof.

For case i), $V$ is continuously differentiable at $x$ by (\ref{eq:DVjsm}) and therefore $DV(x;\eta)=\langle \nabla V(x), \eta \rangle$ is affine with respect to $\eta$. As a result, $DV(x;\dot{x})=\sum_{i \in \I_{sm}(x)} \alpha_i(x) DV(x;f_i(x))$. By the definition of $\I_{sm}$ in (\ref{eq:Ism}) and the construction of $\nu$ in (\ref{eq:Dij}), we have $DV(x;f_i(x))=\min_{q \in \Q} DV(x;f_q(x)), \forall i \in \mathcal{I}_{sm}(x)$. Then the desired condition~(\ref{eq:psclfcond}) follows directly from the decreasing condition of $V$ in (\ref{cond:inf}).

For case ii), According to condition (\ref{eq:bdcond}), there exists some $j^* \in J(x)$ such that 
\begin{align*}
DV_{j^*}(x;f_q(x))\le -W(x), \forall q \in \M_{j^*}^c(x) \cap \I_{sm}^a(x).
\end{align*}
By part ii) of Lemma \ref{lem:DVj}, 
\begin{align*}
DV_{j^*}(x;f_q(x))\le -W(x), \forall q \in \M_{j^*}(x).
\end{align*}
The above two inequalities give
\begin{align}
DV_{j^*}(x;f_i(x))\le -W(x), \forall i \in \I_{sm}^a(x). \label{eq:dvjstar}
\end{align}
To stay on the nonsmooth boundary $\partial \Omega$, the velocity satisfies $\dot{x} \in \bigcap_{j \in J(x)} \T_x(\partial \Omega_j)$. By Lemma \ref{lem:ddBoundary}, $DV(x;\dot{x})=DV_{j^*}(x;\dot{x})$. As we know, $DV_{j^*}(x;\eta)$ is affine with respect to $\eta$, which gives 
\begin{multline*}
DV(x;\dot{x})=\sum_{i \in \I_{sm}(x)} \alpha_i(x) DV_{j^*}(x;f_i(x)) \\ =\sum_{i \in \I_{sm}^a(x)} \alpha_i(x) DV_{j^*}(x;f_i(x))\le -W(x),
\end{multline*}
where the second equality is due to the fact that $\alpha_i(x)=0$, $\forall i \in \I_{sm}(x)\backslash \I_{sm}^a(x)$ and the last inequality follows from (\ref{eq:dvjstar}). Therefore, we complete the proof for case ii). \qed
\end{pf}

According to Theorem \ref{thm:pwsclf2}, the switching law $\nu = \S[V]$ can guarantee closed-loop stability including sliding motions if $V$ is a PSCLF and condition~(\ref{eq:bdcond}) is satisfied on $\partial\Omega\cap\partial\D$. For any given PSCLF, the switching law can be constructed accordingly and condition (\ref{eq:bdcond}) can then be checked without requiring further runtime trajectory dependent information. Note that the index set $\M_j^c(x) \cap \I_{sm}^a(x)$ is a strict subset of $\Q$ and thus the condition~(\ref{eq:bdcond}) is less conservative than checking the inequality for all the subsystems $q\in\Q$. In fact, for many cases, the index set $\M^c_j(x)$ is empty for some $j\in J(x)$, for which (\ref{eq:bdcond}) holds trivially. 

We end this subsection by revisiting Example \ref{ex:DV} to illustrate how to check condition~(\ref{eq:bdcond}) when the switching surface $\partial\D$ partly coincides with the nonsmooth surface $\partial\Omega$. Recall that the nonsmooth surface of $V$ is $S_1\cup S_2$ where $S_1=\{x\in\R^2:x_1-x_2=0\}$ and $S_2=\{x\in\R^2:x_1+3x_2=0\}$. Since $\I_{sm}(x)=\{1,2\}, \forall x\in S_1$, $S_1$ is also a switching boundary, i.e., $S_1 \subseteq \partial \D$. It can be easily verified that $S_2$ is not a switching boundary and therefore $\partial\Omega \cap \partial\D=S_1$. By the relations summarized in (\ref{eq:exdvsmy}), we have $\M_1(x)=\{1\}, \M_2(x)=\{2\}, \forall x\in S_1$. If we fix $j=2$, then the index set to be checked is $\M_2^c(x) \cap \I_{sm}^a(x)=\{1\}$. By the directional derivatives computed in (\ref{eq:exdv}), we have $DV_2(x;f_1(x))=-6\|x\|^2$ and from which condition (\ref{eq:bdcond}) is verified. Therefore, the system is globally asymptotically stable including sliding motions under the switching law $\nu=\S[V]$.

\subsection{Important Special Cases}
This section discusses two important special cases of Theorem~\ref{thm:pwsclf2}. The first one is when the CLF $V$ is smooth over the entire state space, and the second one is when $V$ is obtained by taking the pointwise minimum of a finite number of smooth functions. Both cases have been studied in the literature. Using Theorem~\ref{thm:pwsclf2} we are able to obtain stronger results in a unified way. 

\begin{cor} (Smooth CLF)
\label{cor:smooth}
If $V$ is smooth and satisfies all the conditions (\ref{cond:pd}), (\ref{cond:prop}), and (\ref{cond:inf}), then the closed-loop system~(\ref{eq:cl}) under the switching law $\nu=\S[V]$ is globally asymptotically stable including sliding motions. 
\end{cor}

When $V$ is smooth, the nonsmooth boundary $\partial\Omega$ is empty and condition~(\ref{eq:bdcond}) holds trivially. Thus, the above corollary follows immediately from Theorem~\ref{thm:pwsclf2}. It is worth mentioning that with smooth $V$, the closed-loop vector field under $\nu = \S[V]$ can still be discontinuous with trajectories involving sliding motions. Therefore, Corollary~\ref{cor:smooth} is not a direct consequence of the classical CLF results. 

We next consider a special class of PSCLFs that are obtained by taking the pointwise minimum over a finite number of smooth functions. 

\begin{defn} (Pointwise Minimum CLF) Let $\Vpj:\R^n \to \R_+$, $j\in\N_m$, be nonnegative continuously differentiable functions on $\R^n$. The function defined by 
\begin{align}
\Vp(x)\triangleq \min_{j\in\N_m} \Vpj(x), \quad x\in\R^n, \label{eq:vpm}
\end{align}
is called a pointwise minimum control-Lyapunov function (PMCLF) if
\begin{enumerate}
\item $\Omega_{\text{pm},j}\neq\emptyset$, $\forall j\in\N_m$ where $\Omega_{\text{pm},j}\triangleq \{x\in\R^n: \Vpj(x)<V_{\text{pm},k}(x), \forall k\neq j\}$;
\item $\Vp$ satisfies conditions~(\ref{cond:pd}), (\ref{cond:prop}), and (\ref{cond:inf}).
\end{enumerate}
\label{def:pwminf}
\end{defn}

The first condition in the above definition ensures that every smooth function $\Vpj$ contributes nontrivially to the pointwise minimum. Due to the smoothness of $\Vpj$, each set $\Omega_{\text{pm},j}$ is open with continuously differentiable boundary $\partial\Omega_{\text{pm},j}$. Therefore, a PMCLF $\Vp$ is always piecewise smooth and it is also a PSCLF. 

\begin{cor} (PMCLF Theorem) If $\Vp$ is a PMCLF, then the closed-loop system~(\ref{eq:cl}) under the switching law $\nu = \S[\Vp]$ is globally asymptotically stable including sliding motions. 
\label{cor:pmclf}
\end{cor}

\begin{pf}
Since $\Vp$ is piecewise smooth, it suffices to show that $\Vp$ satisfies condition~(\ref{eq:bdcond}). We pick an arbitrary point $x \in \partial \Omega \cap \partial \mathcal{D}$. Since $\I_{sm}(x)\neq\emptyset$, part ii) of Lemma~\ref{lem:Mj} guarantees the existence of a $j\in J(x)$ such that $\M_j(x) \cap \I_{sm}(x)\neq\emptyset$. We now fix this $j$ and pick an arbitrary $q\in \M_j(x) \cap \I_{sm}(x)$. If $\M_j^c(x)\cap\I_{sm}(x)=\emptyset$, then condition~(\ref{eq:bdcond}) holds trivially. Now assume that $\M_j^c(x)\cap\I_{sm}(x)\neq\emptyset$ and pick an arbitrary $q^\prime\in \M_j^c(x)\cap\I_{sm}(x)$. Then, it follows again from part ii) of Lemma~\ref{lem:Mj} that there exists a $j^\prime \in J(x), j' \neq j$ such that $q^\prime\in \M_{j^\prime}(x) \cap \I_{sm}(x)$. Next, we want to show that for $x\in\partial \Omega$, if $i \in \M_k(x)\cap\I^a_{sm}(x)$ for some $k \in J(x)$, then
\begin{align}
DV_{\text{pm},k} (x;f_i(x)) \ge D\Vpj (x;f_i(x)), \forall j \in J(x). \label{eq:DVpj}
\end{align}
By the regularity condition in Definition \ref{def:regular}, there exists $\epsilon_i>0$ such that $\forall \delta \in (0,\epsilon_i)$, $x-\delta f_i(x) \in \D_i$. It follows from $x \in \partial \Omega_k \cap \partial \D_i$ that $x-\delta f_i(x) \in \Omega_k$. By the definition of $\Vp$ in (\ref{eq:vpm}), $V_{\text{pm},k} (x-\delta f_i(x))\le \Vpj (x-\delta f_i(x)), \forall j \in \N_m$. According to the continuity of $\Vp$, $V_{\text{pm},k} (x)=\Vpj (x)$, $\forall j \in J(x)$. It follows that $\lim_{\delta \downarrow 0} \frac{1}{\delta} (V_{\text{pm},k} (x)-V_{\text{pm},k} (x-\delta f_i(x)))\ge \lim_{\delta \downarrow 0} \frac{1}{\delta} (\Vpj(x)-\Vpj(x-\delta f_i(x)))$ and from which we proved (\ref{eq:DVpj}). Apply (\ref{eq:DVpj}) to some $q' \in \M_j^c(x) \cap \I_{sm}^a(x)$ where $q' \in \M_{j'}(x)$, we have
\begin{align}
D\Vpjj(x;f_{q'}(x)) \ge D\Vpj(x;f_{q'}(x)). \label{eq:DVpjprime}
\end{align}
From part ii) of Lemma \ref{lem:DVj}, $D\Vpjj(x;f_{q'}(x)) \le -W(x)$. Together with (\ref{eq:DVpjprime}), $D\Vpj(x;f_{q'}(x)) \le -W(x)$, $\forall q' \in \M_j^c(x) \cap \I_{sm}^a(x)$, which verifies condition (\ref{eq:bdcond}) and therefore completes the proof. \qed
\end{pf}
Corollary~\ref{cor:pmclf} indicates that if a PMCLF is used, closed-loop stability including sliding motions can always be guaranteed without any extra condition on $\partial\D\cap\partial\Omega$. Such a result holds for any switched nonlinear system and any PMCLF (not necessarily piecewise quadratic). It represents an important contribution on its own.  

\section{Application Examples}
The PSCLF approach, along with the stability results, provides a unified framework to design stabilizing switching laws with a systematic consideration of sliding motions. Once a PSCLF is found, the design of the switching law and the stability analysis of the closed-loop system follow directly from our PSCLF results. The search for a PSCLF can often be done numerically through proper parametrization of the PSCLF. Similar ideas have been studied extensively for switched linear systems (SLSs) with quadratic or piecewise quadratic CLFs~\cite{LM99, HML08}. In this section, we will first briefly show that the proposed framework can be used to recover and extend many existing methods for SLSs in a unified way, and then we will use a numerical example to illustrate its application in stabilization of switched nonlinear systems. 

\subsection{Applications in Switched Linear Systems}
We first consider a general switched linear system (SLS) given by:
\begin{align}
\dot{x}(t)=A_{\sigma(t)}x(t), \sigma(t) \in \mathcal{Q}=\{1, \cdots, M\}, t \in \R_+,
\label{eq:sls}
\end{align}
where $\sigma:\mathbb{R}_+ \to \mathcal{Q}$ denotes the switching signal and $\{A_i\}_{i \in \mathcal{Q}}$ are constant matrices. Note that asymptotic stability is equivalent to exponential stability for SLSs \cite{liberzon2003switching}. This fact is implicitly used in some parts of the following discussions. 

\subsubsection{Quadratic Switching Stabilization}
A well studied stabilization problem for SLSs is the so-called quadratic stabilization problem. System~(\ref{eq:sls}) is called {\em quadratically stabilizable} if there exists a switching law $\nu$ under which the closed-loop system has a quadratic Lyapunov function~\cite{LA09, SESP99}. Using our framework, we call a SLS quadratically stabilizable if it admits a quadratic CLF of the form $V(x) = x^TPx$, $x\in\R^n$. In this special case, it can be easily verified that the PSCLF conditions in (\ref{cond:pd}), (\ref{cond:prop}), and (\ref{cond:inf}) are equivalent to 
\begin{align}
P\succ 0, \text{ and } \min_{i\in\Q} x^T(A_i^TP + PA_i)x <0, \forall x\neq 0.
\label{eq:qs}
\end{align}
This coincides with the strict completeness condition proposed in~\cite{SESP99}. It can be easily verified that the following condition is a sufficient condition to ensure~(\ref{eq:qs}):
\begin{multline}
P\succ 0, \text{ and } \Big(\sum_{i\in\Q}\alpha_i A_i\Big)^TP + P\Big(\sum_{i\in\Q}\alpha_i A_i\Big)\prec 0, \\ \text{for some } \alpha_i\ge 0 \text{ with } \sum_{i\in\Q}\alpha_i = 1.
\end{multline}
The above condition is known as the ``stable convex combination'' condition \cite{LM99, WPD98} and can be used to find $P$ by solving a linear matrix inequality (LMI) feasibility problem. Once a quadratic function satisfying~(\ref{eq:qs}) is found, the stabilizing switching law can be obtained immediately using our framework. Since the function $V(x) = x^TPx$ is globally smooth, the switching law design reduces to $\nu(x) = \argmin_{i\in\Q} x^T(A_i^TP + PA_i)x$. Corollary~\ref{cor:smooth} guarantees the closed-loop system under this $\nu$ is asymptotically stable {\em including sliding motions}. Therefore, our framework can be used to recover most of the existing results in quadratic switching stabilization problems.

\subsubsection{Piecewise Quadratic Switching Stabilization}
As a natural extension of quadratic stabilization, we can consider piecewise quadratic functions as candidate CLFs. A well known result along this direction is the {\em largest-region switching strategy}~\cite{LA09, P03}, whose construction depends on two key components. The first one is a collection of regions defined by $\Omega_i = \{x\in\R^n: x^T H_i x\ge 0\}, i\in\Q$. The second component is a collection of quadratic Lyapunov-like functions $V_i(x) =x^TP_ix$, $x\in\R^n$, $i\in\Q$. Note that for each $i\in\Q$, the matrices $H_i$ and $P_i$ are symmetric but may not be positive or negative semidefinite. Given the two components, the largest region switching strategy is defined by $\nu(x) = \argmax_{i\in\Q} x^TH_ix$. It has been shown that this switching law guarantees closed-loop stability {\em excluding sliding motions} under the following four conditions.
\begin{enumerate}[label=(Q\arabic*')]
\item The union of the regions covers the entire space, i.e., $\cup_{i\in\Q} \Omega_i=\R^n$; 
\item $V_i$ is positive definite on $\Omega_i$; 
\item The Lie derivative of $V_i$ along the vector field of subsystem $i$ is negative definite on $\Omega_i$;
\item  $V_i(x) = V_j(x)$ on $\{x\in\R^n: x^TH_ix = x^TH_jx\}$, for all $i, j\in\Q$.
\end{enumerate}
Using the $\S$-procedure \cite{boyd1994linear}, the matrices $\{H_i, P_i\}_{i\in\Q}$ and hence the largest-region switching strategy can be found by solving some bilinear matrix inequalities (BMIs). The derivation of these BMIs can be found in~\cite{LA09,P03}. 

The largest-region switching strategy mentioned above can be quite conservative. First, the number of Lyapunov-like functions has to be equal to the number of subsystems, which is a restrictive assumption. Second, the decreasing condition for $V_i$ in region $\Omega_i$ is also conservative. Note that the regions $\{\Omega_i\}_{i\in\Q}$ are not mutually exclusive and they differ from the actual switching regions $\hat\Omega_i\triangleq\{x\in\R^n: x^TH_ix>x^TH_jx, \forall j\neq i\}, i\in\Q$. Therefore, requiring $V_i$ to decrease on $\hat\Omega_i$ would be a better choice. Last, selecting subsystem based on the region matrices, although can simplify stability analysis, is less effective. Roughly speaking, the switching control should be chosen to decrease Lyapunov-like functions to achieve a better closed-loop stability performance.

Our framework can be used to tackle these issues and extend the largest-region method in a systematic way. For example, we can consider a piecewise quadratic function $V:\R^n\to\R_+$ with $m$ partitions defined by $\hat\Omega_j=\{x\in\R^n: x^TH_jx>x^TH_kx, \forall k\neq j, k\in\N_m\}, j\in \N_m$. The restriction of $V$ to $\hat\Omega_j$ is assumed to take a quadratic form $V_j(x)=x^TP_jx$, $x\in\hat\Omega_j$. It can be verified that this function $V$ will be a PSCLF for system~(\ref{eq:sls}) if the following conditions hold.
\begin{enumerate}[label=(Q\arabic*)]
\item \label{cond:pwquad1} For each $j\in\N_m$, $V_j$ is positive definite on $\hat\Omega_j$; 
\item \label{cond:pwquad2} For each $j\in\N_m$, $\min_{i\in\Q} x^T(A_i^TP_j + P_jA_i)x< 0, \forall x\in\hat\Omega_j$; 
\item \label{cond:pwquad3} $V_k(x) =V_j(x)$, $\forall x \in\partial\hat\Omega_k\cap\partial\hat\Omega_j$, for all $k, j\in\N_m$. 
\end{enumerate}
When $V$ is a PSCLF, Theorem~\ref{thm:pwsclf1} guarantees that the closed-loop system under the switching law $\nu =\S[V]$ is asymptotically stable excluding sliding motions. Stable sliding motion can be also guaranteed by an additional condition 
\begin{enumerate}[resume,label=(Q\arabic*)]
\item \label{cond:pwquad4} $\forall j\in\N_m, \exists k\in\N_m$ such that $x^T(A_i^TP_k + P_kA_i)x<0,\forall x\in\partial\hat\Omega_j,\forall i\in\Q$.
\end{enumerate}
The above condition \ref{cond:pwquad4} guarantees~(\ref{eq:bdcond}) for the PSCLF $V$ as required by Theorem~\ref{thm:pwsclf2}. A set of BMIs can be derived using the $\S$-procedure to guarantee these four conditions.
\begin{thm} \label{thm:pwquad} If there exists real, symmetric matrices $P_j, H_j, j\in\N_m$, and real numbers $\xi_{jk},\lambda_{ijkt}\in\R$, $\eta_1,\eta_2>0$, $\gamma_{jk},\beta_{jk},\zeta_{ijk}\ge 0$, $\alpha_{ij}\in[0,1], j,k,t\in\N_m,i\in\Q$ such that $\sum_{i\in\Q}\alpha_{ij}=1,\forall j\in\N_m$, and
\begin{subequations}
\begin{align}
& P_j-\eta_1 I \succeq \sum_{k\in\N_m} \gamma_{jk}(H_j-H_k), \forall j\in\N_m; \label{eq:pwquad1} \\
& \begin{aligned}\Big( \sum_{i\in\Q} \alpha_{ij} A_i \Big)^TP_j + P_j\Big( \sum_{i\in\Q} \alpha_{ij} A_i \Big)+\eta_2P_j \\ \preceq \sum_{k\in\N_m} \beta_{jk}(H_k-H_j),\forall j\in\N_m; \end{aligned} \label{eq:pwquad2} \\
& P_j=P_k+\xi_{jk}(H_j-H_k),\forall j,k\in\N_m; \label{eq:pwquad3} \\
& \begin{aligned} & \sum_{k\in\N_m} \zeta_{ijk}(A_i^TP_k+P_kA_i) + \\ & \sum_{k,t\in\N_m} \lambda_{ijkt}(H_j-H_t) \prec 0,\forall i\in\Q, \forall j\in\N_m, \end{aligned} \label{eq:pwquad4}
\end{align}
\end{subequations}
then under the switching law $\nu=\S[V]$, system (\ref{eq:sls}) is globally asymptotically stable \emph{including sliding motions}.
\end{thm}
\begin{pf}
By the $\S$-procedure, it is easy to verify that (\ref{eq:pwquad1}) and (\ref{eq:pwquad3}) imply condition \ref{cond:pwquad1} and \ref{cond:pwquad3}, respectively. For any $x\in\hat{\Omega}_j$, we have $x^TH_jx\ge x^TH_kx, \forall k\in\N_m$ and hence $\beta_{jk}x^T(H_k-H_j)x\le 0$ for any $\beta_{jk}\ge 0$. By left multiplying $x^T$ and right multiplying $x$ on both sides of (\ref{eq:pwquad2}), we have
\begin{align*}
DV_j\Big(x;\sum_{i\in\Q} \alpha_{ij} A_ix \Big)+\eta_2 V_j(x) \le 0, \forall x\in\hat{\Omega}_j.
\end{align*}
Since $\sum_{i\in\Q}\alpha_{ij}=1$, it follows that for each $j\in\N_m$,
\begin{multline*}
\min_{i\in\Q} DV_j(x;A_ix)\le DV_j\Big(x;\sum_{i\in\Q} \alpha_{ij} A_ix \Big) \\ \le -\eta_2V_j(x),\forall x\in\hat{\Omega}_j,
\end{multline*}
which verifies condition \ref{cond:pwquad2}. Next, we want to show that (\ref{eq:pwquad4}) implies condition \ref{cond:pwquad4}. Since $\zeta_{ijk}\ge 0$,
\begin{align*}
\sum_{k\in\N_m} \zeta_{ijk}(A_i^TP_k+P_kA_i - \sum_{t\in\N_m} b_{ijkt}(H_j-H_t)) \prec 0
\end{align*}
implies that there exists at least one $k\in\N_m$ such that $A_i^TP_k+P_kA_i - \sum_{t\in\N_m} b_{ijkt}(H_j-H_t) \prec 0$, otherwise the above inequality changes direction. We now fix this $k$ and by the $\S$-procedure we know that $x^T(A_i^TP_k+P_kA_i)x<0,\forall x\in \partial\hat{\Omega}_j$, which verifies condition \ref{cond:pwquad4}. Let $\lambda_{ijkt}\triangleq -\zeta_{ijk}b_{ijkt}$, we obtain the form in (\ref{eq:pwquad4}). \qed
\end{pf}

The search for a piecewise quadratic CLF can thus be formulated as a feasibility problem for the BMIs given in Theorem~\ref{thm:pwquad}. Such a result is more general than the largest-region switching approach~\cite{P03} as it allows the number of the switching regions to be different from the number of subsystems, and guarantees closed-loop stability including sliding motions. 

\subsubsection{Switching Stabilization with Composite Control-Lyapunov Functions}
Another class of CLFs that have been studied for SLSs is composite quadratic functions, which are defined by taking the pointwise minimum, pointwise maximum, and convex hull of a finite number of quadratic functions~\cite{HML08}. The stabilization results based on these three classes of composite CLFs are also special cases of our PSCLF framework. Moreover, they can be further extended and strengthened using our framework. As an example, we consider a pointwise maximum CLF defined by:
\begin{align*}
\Vpm\triangleq \max_{j\in\N_m} V_j(x), \quad x\in\R^n,
\end{align*}
where $V_j(x)=x^TP_jx$. Obviously, $\Vpm$ is a piecewise quadratic function with $m$ partitions defined by $\Omega_j=\{x\in\R^n: x^TP_jx>x^TP_kx, \forall k\neq j, k\in\N_m\}, j\in \N_m$. It can be verified that $\Vpm$ will be a PSCLF for system (\ref{eq:sls}) if the following conditions hold.
\begin{enumerate}[label=(M\arabic*)]
\item \label{cond:ptmax1} For each $j\in\N_m$, $V_j$ is positive definite on $\Omega_j$; 
\item \label{cond:ptmax2} For each $j\in\N_m$, $\min_{i\in\Q} x^T(A_i^TP_j + P_jA_i)x< 0, \forall x\in\Omega_j$. 
\end{enumerate}
When $\Vpm$ is a PSCLF, Theorem \ref{thm:pwsclf1} guarantees that the closed-loop system under the switching law $\nu=\S[\Vpm]$ is asymptotically stable excluding sliding motions. Stable sliding motion can be guaranteed by an additional condition
\begin{enumerate}[resume,label=(M\arabic*)]
\item \label{cond:ptmax3} $\forall j\in\N_m, \exists k\in\N_m$ such that $x^T(A_i^TP_k + P_kA_i)x<0,\forall x\in\partial\Omega_j,\forall i\in\Q$.
\end{enumerate}
Note that the above condition implies condition (\ref{eq:bdcond}). It requires $x\in\partial \Omega_j$ rather than $x\in\partial\Omega_j\cap\partial\D_i$ for ease of BMI derivation. The construction of $\Vpm$ and the corresponding switching law design can be formulated as a BMI feasibility problem as stated in the following theorem.
\begin{thm} If there exists real, symmetric matrices $P_j, j\in\N_m$, and real numbers $\lambda_{ijkt}\in\R$, $\eta_1,\eta_2>0$, $\gamma_{jk},\beta_{jk},\zeta_{ijk}\ge 0$, $\alpha_{ij}\in[0,1], j,k,t\in\N_m,i\in\Q$ such that $\sum_{i\in\Q}\alpha_{ij}=1,\forall j\in\N_m$, and
\begin{subequations}
\begin{align}
& P_j-\eta_1 I \succeq \sum_{k\in\N_m}\gamma_{jk}(P_j-P_k),\forall j\in\N_m; \label{eq:ptmax1}\\
& \begin{aligned}\Big( \sum_{i\in\Q} \alpha_{ij} A_i \Big)^TP_j + P_j\Big( \sum_{i\in\Q} \alpha_{ij} A_i \Big)+\eta_2P_j \\ \preceq \sum_{k\in\N_m} \beta_{jk}(P_k-P_j),\forall j\in\N_m; \end{aligned} \label{eq:ptmax2} \\
& \begin{aligned} & \sum_{k\in\N_m}\zeta_{ijk} (A_i^TP_k+P_kA_i)+ \\ & \sum_{k,t\in\N_m} \lambda_{ijkt}(P_j-P_t) \prec 0,\forall i\in\Q, \forall j\in\N_m, \end{aligned} \label{eq:ptmax3}
\end{align}
\end{subequations}
then under the switching law $\nu=\S[\Vpm]$, system (\ref{eq:sls}) is globally asymptotically stable \emph{including sliding motions}.
\end{thm}
\begin{pf}
By the $\S$-procedure, condition \ref{cond:ptmax1} is guaranteed by (\ref{eq:ptmax1}). Similar as the argument of (\ref{eq:pwquad2}) $\Rightarrow$ condition \ref{cond:pwquad2} shown in the proof of Theorem \ref{thm:pwquad}, we have (\ref{eq:ptmax2}) implies condition \ref{cond:ptmax2}. We are left to show that (\ref{eq:ptmax3}) implies condition \ref{cond:ptmax3}. Since $\zeta_{ijk}\ge 0$, $\sum_{k\in\N_m} \zeta_{ijk}(A_i^TP_k+P_kA_i - \sum_{t\in\N_m}b_{ijkt}(P_j-P_t)) \prec 0$ implies that there exists at least one $k\in\N_m$ such that
\begin{align}
A_i^TP_k+P_kA_i - \sum_{t\in\N_m}b_{ijkt}(P_j-P_t) \prec 0. \label{eq:ptmaxbd}
\end{align}
We now fix this $k$. Apply the $\S$-procedure to (\ref{eq:ptmaxbd}), we have $x^T(A_i^TP_k + P_kA_i)x< 0$ on $\{x\in\R^n: x^TP_jx=x^TP_tx,t\in\N_m\}=\cup_{t\in\N_m} (\partial\Omega_j\cap\partial\Omega_t)$. Due to the fact that $\partial\Omega_j=\cup_{t\in\N_m} (\partial\Omega_j\cap\partial\Omega_t)$, condition \ref{cond:ptmax3} is verified. Let $\lambda_{ijkt}\triangleq -\zeta_{ijk}b_{ijkt}$, we obtain the form in (\ref{eq:ptmax3}). \qed
\end{pf}

The authors in \cite{HML08} also derived a set of conditions to guarantee closed-loop stability including sliding motions. The conditions are mostly the same as ours except for the one on the nonsmooth boundary $\partial\Omega$. The BMI conditions in \cite{HML08} are derived only for the case where $\Vpm$ is composed from two quadratic functions (i.e. $m=2$), which is more restrictive than condition (\ref{eq:ptmax3}) obtained using our framework. It is worth mentioning that similar extensions can be made for the case with convex-hull composite CLFs. As for the pointwise minimum CLF case, the result in \cite{HML08} is a very special case of our general result given in Corollary~\ref{cor:pmclf}. All these discussions further illustrate the importance and unified nature of the proposed PSCLF framework. 


\begin{figure*}[!ht]
	\centering
	\begin{tabular}{ccc}
		\includegraphics[width=0.3\textwidth]{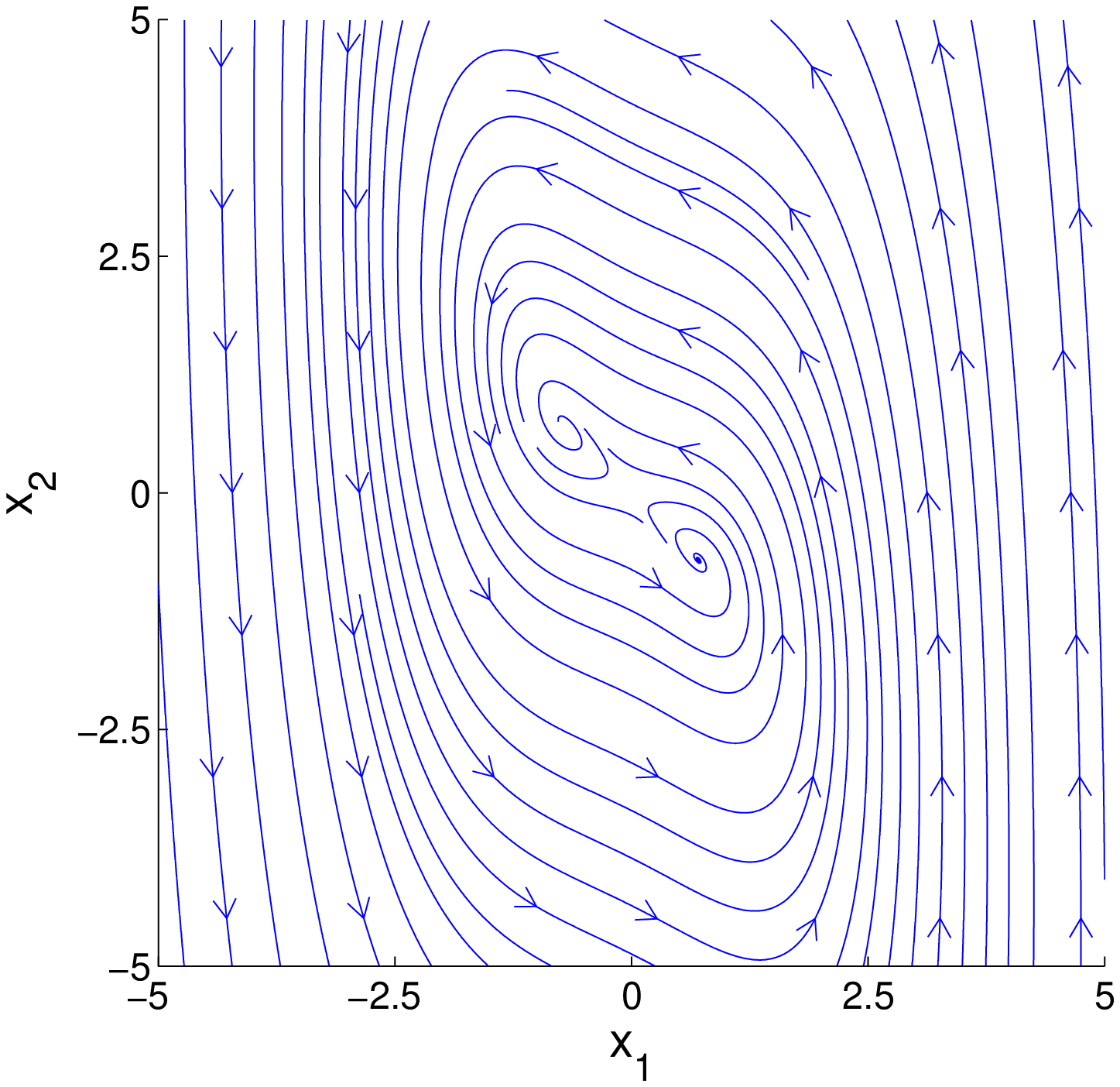} &
		\includegraphics[width=0.3\textwidth]{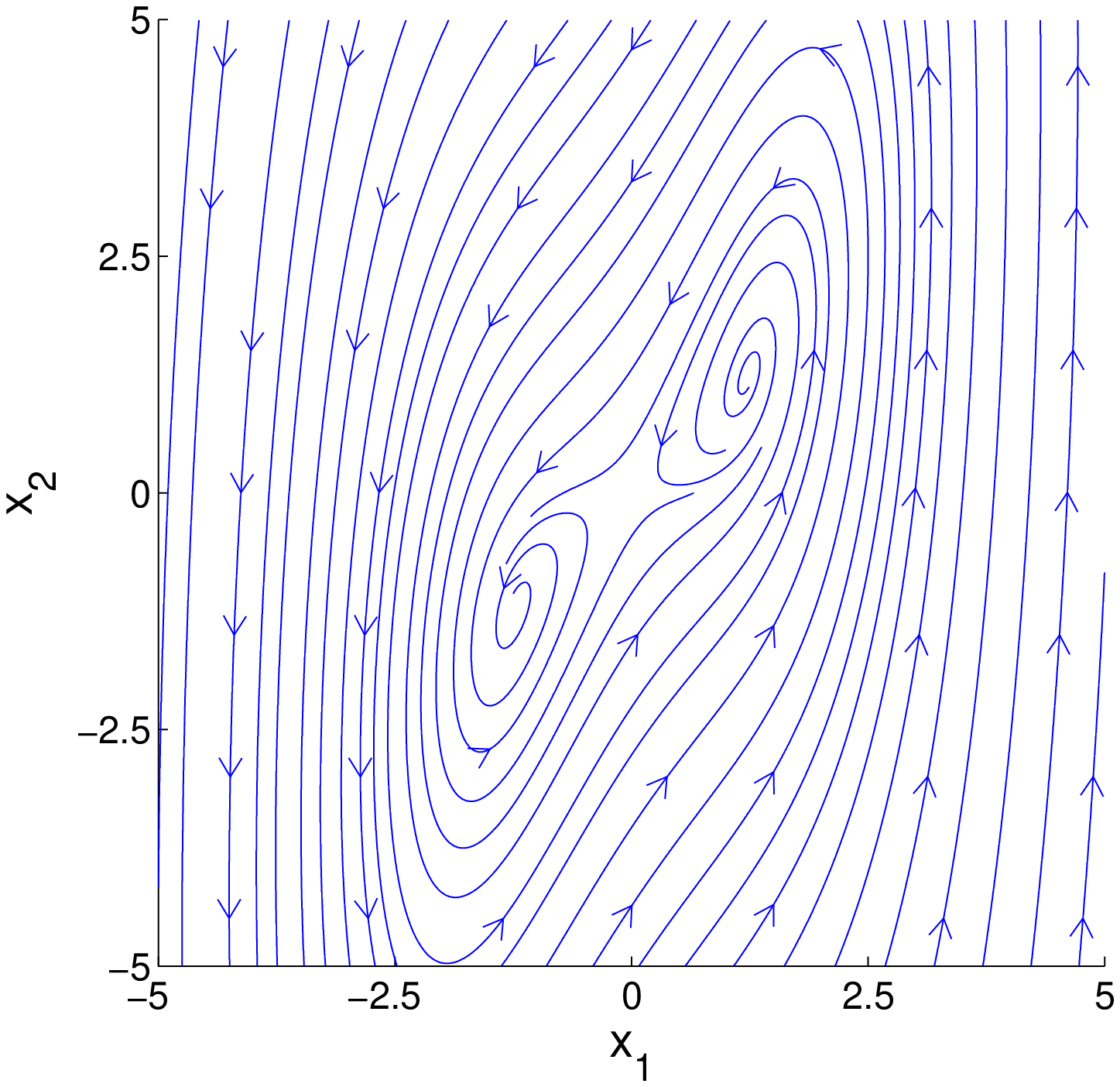} &
		\includegraphics[width=0.3\textwidth]{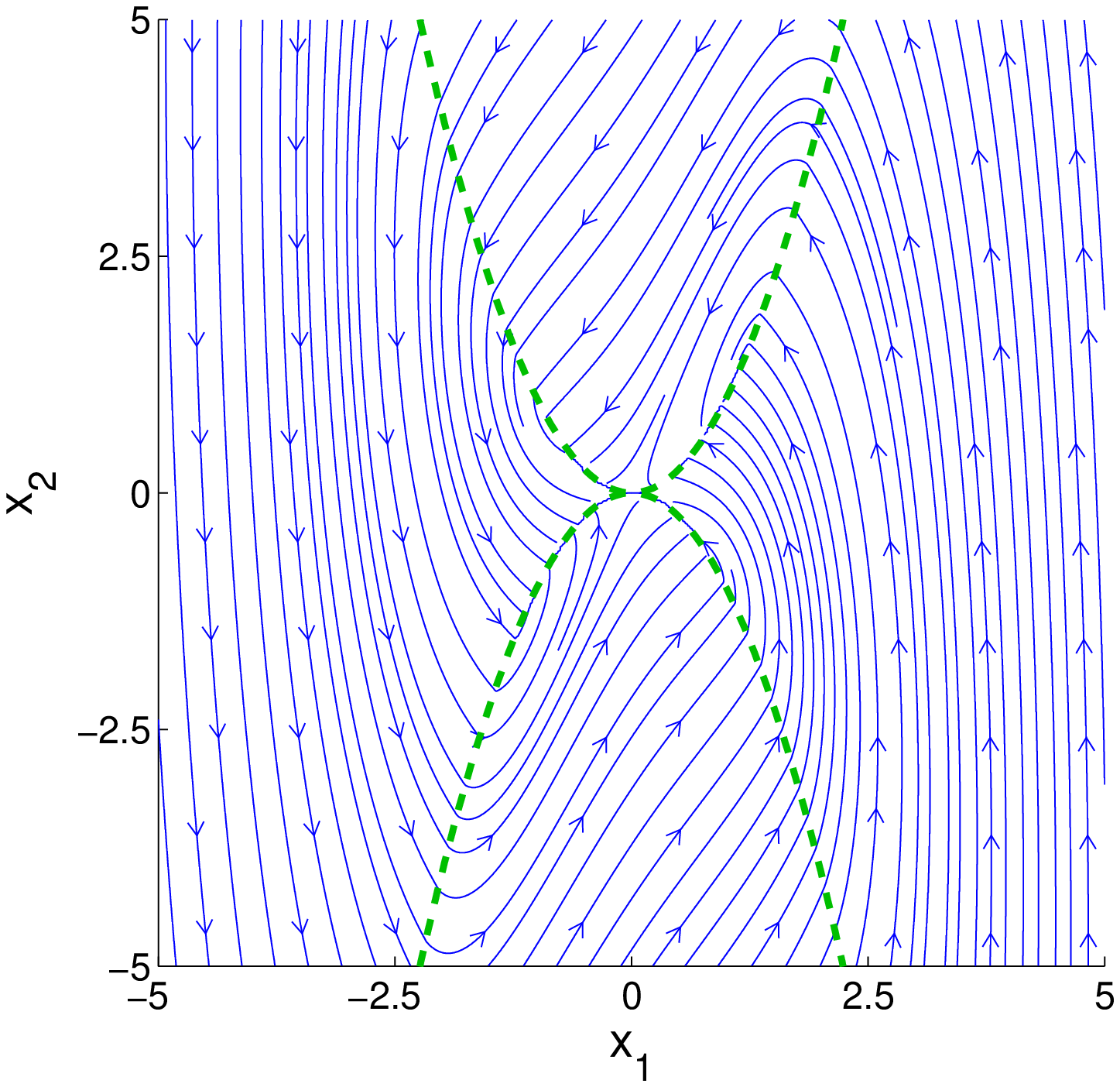} \\
		(a) & (b) & (c)
	\end{tabular}
	\begin{tabular}{cc}
		\includegraphics[width=0.4\textwidth]{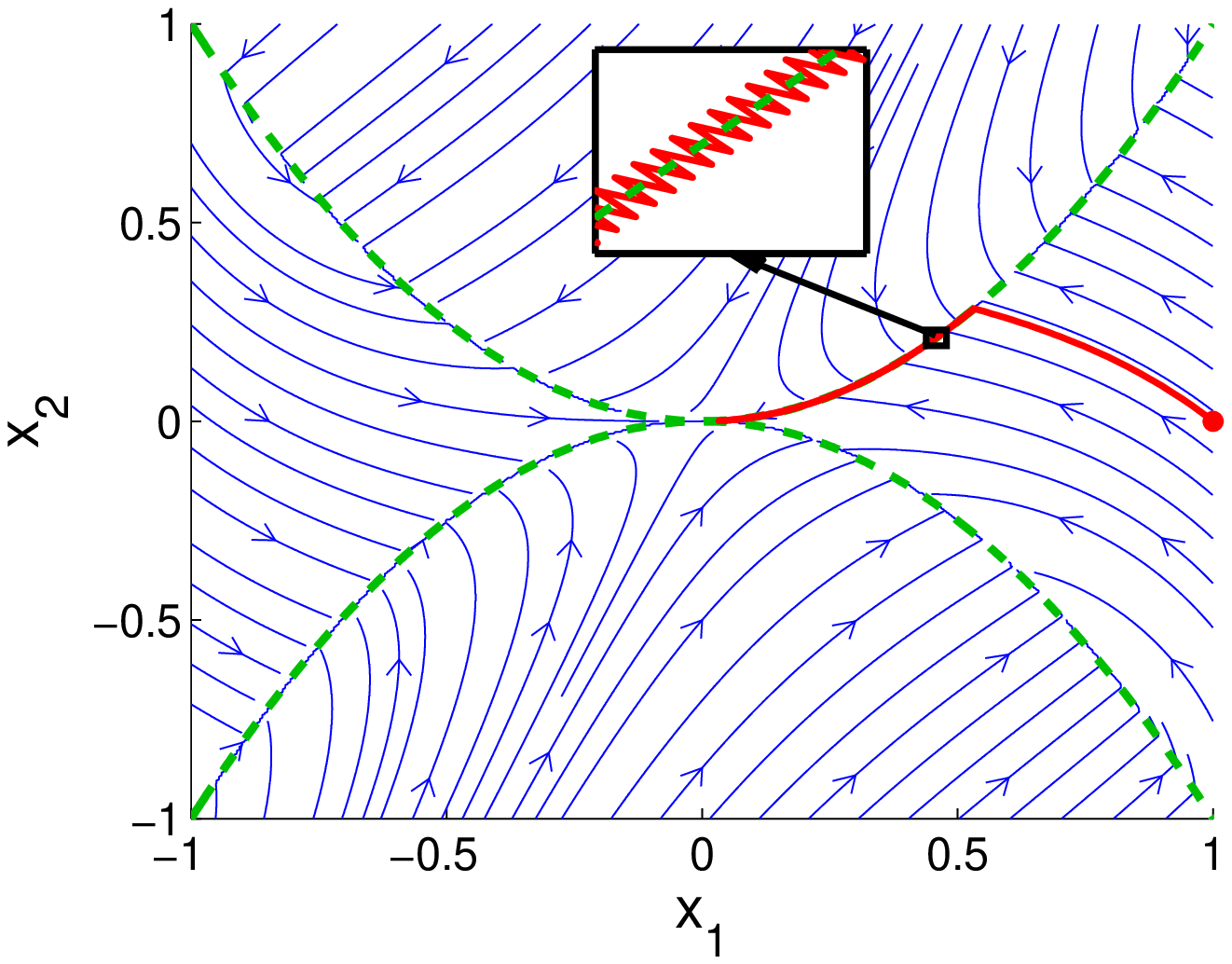} &
		\includegraphics[width=0.4\textwidth]{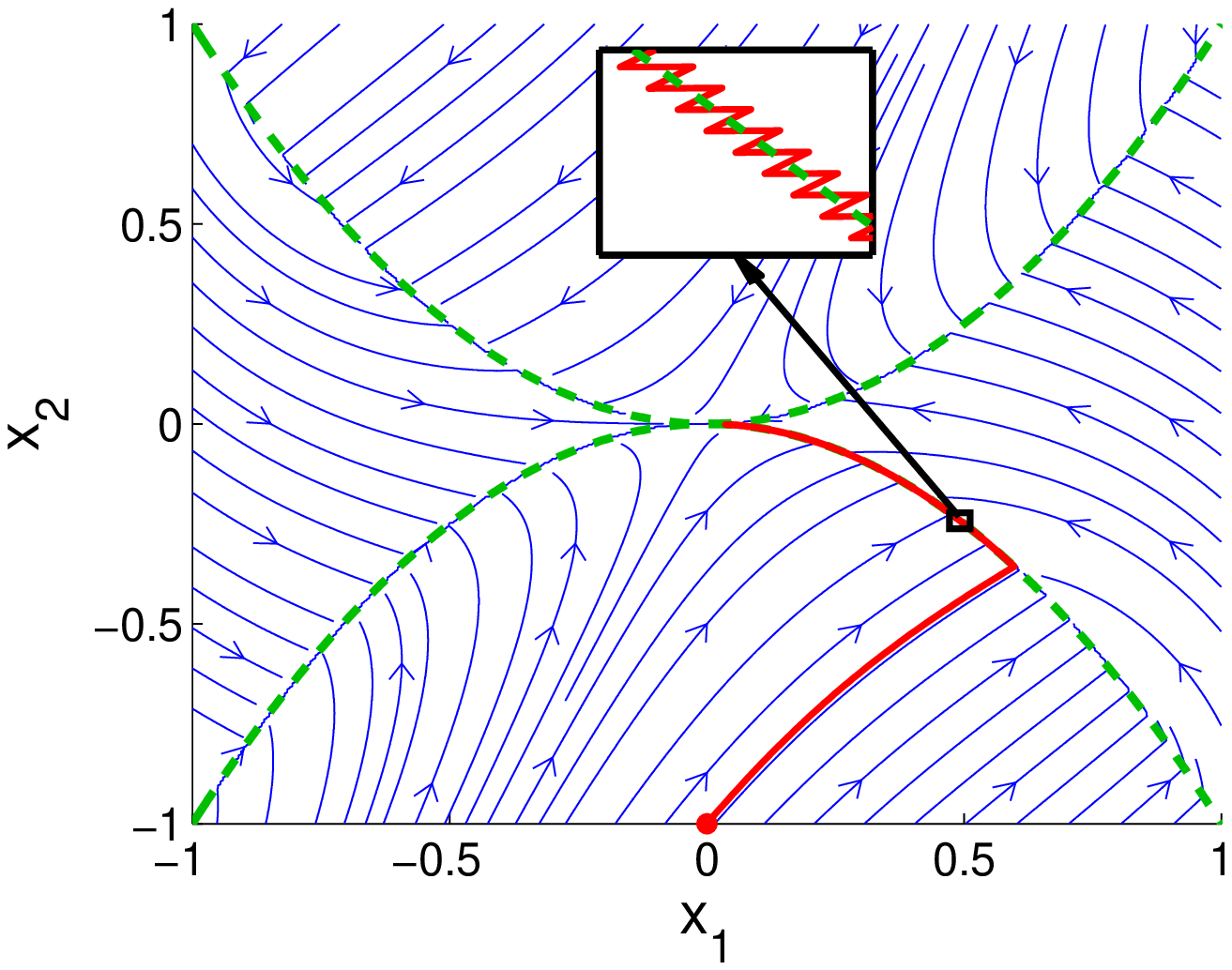} \\
		(d) & (e)	
	\end{tabular}
	\caption{(a): The phase portrait of subsystem 1, (b): The phase portrait of subsystem 2, (c): The phase portrait of the closed-loop system under the switching law $\nu$, (d): The closed-loop trajectory of system (\ref{eq:nonlinstab}) with $x_0=[1,0]^T$, (e): The closed-loop trajectory of system (\ref{eq:nonlinstab}) with $x_0=[0,-1]^T$.}
	\label{fig:nonlin}
\end{figure*}

\subsection{A Numerical Nonlinear Example}
Switching stabilization of switched nonlinear systems has not been adequately studied in the literature. Here, we use a numerical example to illustrate the application of our framework in this area. Consider the following switched nonlinear system with two subsystems, $\dot{x}(t)=f_i(x(t)), i=1,2$, where
\begin{align}
f_1(x)=\begin{bmatrix} -x_1-x_2 \\ x_1^3+0.5x_2 \end{bmatrix}, 
f_2(x)=\begin{bmatrix} x_1-x_2 \\ x_1^3-1.5x_2 \end{bmatrix}.
\label{eq:nonlinstab}
\end{align}

For this simple switched polynomial system with subsystem vector fields shown in Fig. \ref{fig:nonlin}(a) and Fig. \ref{fig:nonlin}(b), we consider a smooth polynomial CLF: $V(x)=x_1^4+2x_2^2$ under which the closed-loop vector field is shown in Fig. \ref{fig:nonlin}(c). Apparently, $V$ satisfies the positive definite condition and the radially unbounded condition given in (\ref{cond:pd}) and (\ref{cond:prop}). In addition, we have $DV(x;f_1(x))=-(4x_1^4-4x_2^2)-2x_2^2$, $DV(x;f_2(x))=(4x_1^4-4x_2^2)-2x_2^2$. The above implies $\min_i DV(x;f_i(x))\le -2x_2^2$, which verifies the decreasing condition in (\ref{cond:inf}). Therefore, $V$ is a special case of the proposed PSCLFs. Since it is smooth, the corresponding switching law in (\ref{eq:nu}) takes the following form: $\nu(x)=\argmin_i DV(x;f_i(x))$. By Corollary \ref{cor:smooth}, we can conclude that the closed-loop system $\dot{x}=f_{\nu(x)}(x)$ is asymptotically stable {\em including sliding motions}. The result is also illustrated through simulations. In particular, Fig.~\ref{fig:nonlin}(d) and Fig.~\ref{fig:nonlin}(e) show the closed-loop trajectories starting from $x_0=[1,0]^T$ and $x_0=[0,-1]^T$, respectively. Both trajectories exhibit sliding motion behavior, which can be observed in the zoom-in box. We introduce hysteresis band in the simulation to deal with the discontinuous closed-loop vector field. Both trajectories converge to the origin under sliding motion as expected.

\section{Conclusion}
In this paper, we proposed a piecewise smooth control-Lyapunov function (PSCLF) framework to study switching stabilization problems. We formally introduced the concept of stability including or excluding sliding motions and sufficient conditions in terms of PSCLFs were derived for the two stability notions. A constructive way to design a stabilizing switching law based on a given PSCLF was also developed. We showed that such a control law can guarantee the closed-loop stability excluding sliding motions, and it can also ensure closed-loop stability including sliding motions under an additional condition on the nonsmooth surface of the PSCLF. We also showed that for smooth CLFs and pointwise minimum CLFs, stable sliding motions are always guaranteed without any additional condition. The proposed framework can be used to obtain and extend most of the existing results in stabilization of switched linear systems. In addition, it provides a systematic way to study stabilization of switched nonlinear systems. Future research will focus on necessary stabilizability conditions and converse control-Lyapunov function theorems for switched linear systems.

\bibliographystyle{abbrv}        
\bibliography{psclf}

\appendix

\end{document}